\documentclass[11pt]{amsart}
\usepackage[margin=1in]{geometry}
\usepackage{amsmath,amssymb,amsthm,mathtools}
\usepackage{enumitem}
\usepackage{booktabs}
\usepackage[colorlinks=true,linkcolor=blue,citecolor=blue]{hyperref}

\newtheorem{theorem}{Theorem}[section]
\newtheorem{lemma}[theorem]{Lemma}
\newtheorem{proposition}[theorem]{Proposition}
\newtheorem{corollary}[theorem]{Corollary}
\newtheorem{mainthm}{Theorem}

\theoremstyle{definition}
\newtheorem{definition}[theorem]{Definition}
\newtheorem*{definition*}{Definition}
\newtheorem{example}[theorem]{Example}
\newtheorem{problem}[theorem]{Problem}
\theoremstyle{remark}
\newtheorem{remark}[theorem]{Remark}

\newcommand{\R}{\mathbb{R}}
\newcommand{\N}{\mathbb{N}}
\newcommand{\eps}{\varepsilon}
\newcommand{\SOC}{\textup{(SOC)}}
\newcommand{\SDP}{\textup{(SDP)}}
\newcommand{\ocv}{\xrightarrow{\,\mathrm o\,}}
\newcommand{\fb}{\mathfrak b}
\newcommand{\OCRU}{\textup{(OC$\Rightarrow$RUC)}}

\begin{document}

\title[Completions and the diagonal property]{Completions of vector metric spaces and the diagonal property of order convergence}

\author{C\"uneyt \c{C}ev\.{\i}k}
\address{Department of Mathematics, Faculty of Science, Gazi University, Ankara, Turkey}
\email{ccevik@gazi.edu.tr}


\begin{abstract}
A vector metric space carries a distance with values in an ordered vector space $E$; a sequence converges if its distances to the limit are dominated by a sequence decreasing to zero in $E$. We show that, for Archimedean $E$, every $E$-valued vector metric space has a completion if and only if order Cauchy sequences in $E$ converge and $E$ has the diagonal property of order convergence \SDP: for every double sequence $(a^k_n)$ with $a^k_n\downarrow_n0$ there are indices $n_k$ and $e_k\downarrow 0$ with $a^k_{n_k}\le e_k$. In vector lattices this is the classical diagonal property of Kantorovich, equivalently the strong Egoroff property. We show that its standard characterisations remain valid in Archimedean ordered vector spaces, and that it holds if and only if convergence in every $E$-valued vector metric space satisfies the iterated limit axiom. For an ordered Banach space with a closed normal cone, all completions exist if and only if the norm is $\sigma$-order continuous; if the cone has an interior point, all vector metrics are then scalar. By a theorem of Veksler and of Hager and van Mill, for compact Hausdorff $K$ either $K$ is an almost $P$-space and $C(K)$-valued vector metrics are scalar, or there are countable $C(K)$-valued vector metric spaces without completion. A product of ordered vector spaces with \SDP\ and nonzero cones has \SDP\ if and only if there are fewer than $\fb$ factors, extending a result of Hager and van Mill for $\R^A$. We also prove permanence results for direct summands and uniformly dense vector lattice covers of pre-Riesz spaces.
\end{abstract}

\maketitle

\section{Introduction}\label{sec:intro}

Distances with values in ordered vector spaces have a long tradition. Kantorovich's $K$-metric spaces are surveyed by Zabrejko \cite{Zab97}; there convergence is defined through a norm or a lattice structure on the value space. \c{C}evik and Altun~\cite{CA09} introduced \emph{vector metric spaces}, in which the distance takes its values in a vector lattice $E$ and convergence is purely order-theoretic: a sequence converges if its distances to the limit are dominated by a sequence decreasing to $0$ in $E$ (Definition~\ref{def:vm}). The definition makes sense verbatim when $E$ is only an ordered vector space. Continuity of maps between such spaces was studied in \cite{Cev14}, and completions and unbounded completions in \cite{Cev15,OC20}; Remark~\ref{rem:Cev15} relates them to the present results.

The construction of a completion by Cauchy sequences meets two obstacles in this setting. First, for Cauchy sequences $(x_n)$ and $(y_n)$ the distances $d(x_n,y_n)$ must converge in $E$; this requires that order Cauchy sequences in $E$ converge, a property we call \SOC. Second, and more seriously, the space of equivalence classes must itself be complete. Classically this is a diagonal argument, and order convergence is notoriously ill-behaved with respect to diagonal arguments. The relevant property of $E$ is the following diagonal property of order convergence.

\begin{definition*}
An ordered vector space $E$ has the \emph{sequential diagonal property} \SDP\ if for every family $(a^k_n)_{k,n\in\N}$ in $E$ with $a^k_n\downarrow_n0$ for each $k$ there are $n_1,n_2,\ldots\in\N$ and a sequence $e_k\downarrow 0$ in $E$ with $a^k_{n_k}\le e_k$ for all $k$.
\end{definition*}

For vector lattices, \SDP\ is a classical notion. It is the diagonal property for order convergence, which goes back to Kantorovich \cite[Satz~24]{Kan37} and is treated in \cite[Chapter~10]{LZ71}. Dodds \cite{Dod74} showed that in a vector lattice it is equivalent to the strong Egoroff property and to the idempotence of the pseudo order closure; for Archimedean vector lattices it is also equivalent to the $\sigma$-property together with stability of order convergence \cite[Theorem~70.2]{LZ71}. In this last form it is the definition of \emph{strongly Egoroff} vector lattices in recent work of Hager and van Mill \cite{HvM15} and Hager \cite{Hag16}, who determined the property for many spaces of continuous functions. For instance, $\R^A$ is strongly Egoroff if and only if $|A|<\fb$, where $\fb$ is the bounding number \cite[Remark~2.2(ii)]{HvM15}, and $C(K)$ for compact $K$ is strongly Egoroff if and only if $K$ is an almost $P$-space \cite[Theorem~2]{Vek73}, \cite[Theorem~6.4]{HvM15}; the weaker Egoroff property of $\R^A$ was determined by Blass and Jech \cite{BJ86}. Remark~\ref{rem:classical} gives a more detailed dictionary.

We keep the abbreviation \SDP, which stresses that only sequences are involved, and we work in ordered vector spaces. Many natural value spaces of vector metrics are not lattices, for instance spaces of smooth functions, spaces of self-adjoint operators and spin factors, and the lattice-theoretic arguments of \cite{LZ71,Dod74,HvM15} do not apply to them. Our results show that the classical property, together with \SOC, is exactly what completions of vector metric spaces require, and that its basic theory does not depend on lattice structure.

We now state the main results. Throughout, $E$ is a real ordered vector space, compact spaces are Hausdorff, and $c_n\downarrow 0$ means that $(c_n)$ is decreasing with infimum $0$.

\begin{mainthm}[Theorem~\ref{thm:main}]\label{thmA}
Let $E\ne\{0\}$ be Archimedean. The following are equivalent:
\begin{enumerate}[label=\textup{(\roman*)}]
\item every $E$-valued vector metric space has a completion;
\item $E$ satisfies \SOC\ and \SDP.
\end{enumerate}
The completion is then unique up to isometry. Moreover, $E$ has \SDP\ if and only if convergence in every $E$-valued vector metric space satisfies the iterated limit axiom.
\end{mainthm}

The necessity part uses explicit countable vector metric spaces (Lemma~\ref{lem:rows}); when \SDP\ fails, some of them have no completion, and some carry a convergence that is induced by no first countable topology. In Section~\ref{sec:SDP} we show that the classical characterisations of \SDP\ by the dominating property and by the $\sigma$-property together with stability of order convergence hold in every Archimedean ordered vector space (Proposition~\ref{prop:egoroff}).

\begin{mainthm}[Theorems~\ref{thm:banach} and \ref{thm:vm}, Proposition~\ref{prop:unit}]\label{thmB}
Let $E$ be an ordered Banach space with a closed normal cone. Then $E$ has \SDP\ if and only if its norm is $\sigma$-order continuous, if and only if every continuous positive linear functional on $E$ is $\sigma$-order continuous; in this case $E$ satisfies \SOC, so completions exist. If moreover the cone has an interior point, then under these conditions order convergence of sequences coincides with norm convergence, and every $E$-valued vector metric space is scalar: its convergence, completeness and completion are those of the metric $\|d\|$.
\end{mainthm}

For Banach lattices the first equivalence also follows from \cite[Theorem~70.2]{LZ71}; the point of Theorem~\ref{thmB} is that no lattice structure is needed. It applies to the non-lattice space $S_2^{\mathrm{sa}}$ of self-adjoint Hilbert--Schmidt operators with the L\"owner order and to all spin factors, and it shows that the self-adjoint part of every infinite-dimensional von Neumann algebra fails \SDP. Outside the Banach setting the functional criterion fails: in $\R^{\fb}$ and in $c_{00}$ all positive functionals are $\sigma$-order continuous, yet \SDP\ fails (Proposition~\ref{prop:converse}). Completions need not exist for values in $\ell^\infty$, although $\ell^\infty$ is Dedekind complete (Example~\ref{ex:linfty}).

\begin{mainthm}[Theorem~\ref{thm:CK} and Corollary~\ref{cor:CKbad}]\label{thmC}
Let $K$ be a compact Hausdorff space.
\begin{enumerate}[label=\textup{(\roman*)}]
\item If $K$ is an almost $P$-space, every $C(K)$-valued vector metric space has a completion, and its convergence and completeness are those of the metric $\|d\|_\infty$.
\item Otherwise there are countable $C(K)$-valued vector metric spaces without completion, and countable ones whose convergence is induced by no first countable topology.
\end{enumerate}
\end{mainthm}

The order-theoretic core of Theorem~\ref{thmC}, that $C(K)$ has \SDP\ exactly for almost $P$-spaces $K$, was announced without proof by Veksler \cite[Theorem~2]{Vek73}, who also lists the diagonal property among the equivalent conditions, and proved by Hager and van Mill \cite[Theorem~6.4]{HvM15}. Thus $C(K)$-valued vector metrics are genuinely vector-valued only in case (ii), where they are badly behaved. Examples for (i) are $C(\N^*)=\ell^\infty/c_0$ and $C(D\cup\{\infty\})$ for uncountable discrete $D$; examples for (ii) are $\ell^\infty=C(\beta\N)$, $c$, and $C(K)$ for every infinite metrisable $K$.

\begin{mainthm}[Theorems~\ref{thm:small} and \ref{thm:large}]\label{thmD}
A product $\prod_{\alpha\in A}E_\alpha$ of ordered vector spaces with \SDP\ and nonzero cones has \SDP\ if and only if $|A|<\fb$.
\end{mainthm}

For $E_\alpha=\R$ this is the result of Hager and van Mill quoted above; our proof is direct and uses no lattice operations. Consequently every $\R^A$-valued vector metric space has a completion if and only if $|A|<\fb$, which for $A=\omega_1$ is independent of ZFC (Corollary~\ref{cor:indep}). In Remark~\ref{rem:CXknown} and Corollary~\ref{cor:CX} we combine the results of \cite{HvM15,Hag16} on spaces $C(X)$ with Theorem~\ref{thmA}.

\begin{mainthm}[Section~\ref{sec:perm}]\label{thmE}
The property \SDP\ passes to order direct summands and from a pre-Riesz space to a vector lattice cover in which it is uniformly dense, and back. It passes neither to nor from order dense ideals, and it need not pass from a vector lattice to its Dedekind completion, even when the latter is a pervasive cover.
\end{mainthm}

The paper is organised as follows. Section~\ref{sec:prelim} collects preliminaries on order convergence of sequences, vector metric spaces and the bounding number. Section~\ref{sec:SDP} contains basic properties of \SDP\ and its relation to the Egoroff properties. Section~\ref{sec:completions} proves Theorem~\ref{thmA}, Section~\ref{sec:banach} treats ordered Banach spaces and spaces $C(K)$, Section~\ref{sec:products} products and spaces $C(X)$, and Section~\ref{sec:perm} permanence properties. Section~\ref{sec:further} contains further examples, a table and open problems.

\section{Preliminaries}\label{sec:prelim}

We write $E_+$ for the cone of $E$, so $E_++E_+\subseteq E_+$, $\lambda E_+\subseteq E_+$ for $\lambda\ge 0$ and $E_+\cap(-E_+)=\{0\}$. The space $E$ is \emph{Archimedean} if $nx\le y$ for all $n\in\N$ implies $x\le 0$. Sequences are indexed by $\N=\{1,2,\dots\}$.

Order convergence of sequences is defined by domination; no suprema or infima of the sequence itself are involved, so the definition makes sense in every ordered vector space.

\begin{definition}\label{def:oconv}
Let $(z_n)$ be a sequence in $E$ and $z\in E$.
\begin{enumerate}[label=(\alph*)]
\item $z_n\ocv z$ (\emph{order convergence}) if there is $c_n\downarrow 0$ with $\pm(z_n-z)\le c_n$ for all $n$.
\item $(z_n)$ is \emph{order Cauchy} if there is $c_n\downarrow 0$ with $\pm(z_n-z_m)\le c_n$ for all $m\ge n$.
\end{enumerate}
\end{definition}

We shall use the following rules throughout; they replace the lattice calculus.

\begin{lemma}\label{lem:olim}
\begin{enumerate}[label=\textup{(\roman*)}]
\item If $c_n\downarrow 0$ and $c'_n\downarrow 0$, then $c_n+c'_n\downarrow 0$, and every lower bound of $\{c_n+c'_m:n,m\in\N\}$ is $\le 0$.
\item Order limits are unique, and $z_n\ocv z$, $w_n\ocv w$ imply $\alpha z_n+\beta w_n\ocv\alpha z+\beta w$ for $\alpha,\beta\in\R$.
\item If $z_n\ocv z$, $w_n\ocv w$ and $z_n\le w_n$ for all large $n$, then $z\le w$.
\item If $\pm(z_n-z_m)\le c_n$ for all $m\ge n$ and $z_n\ocv z$, then $\pm(z_n-z)\le c_n$ for all $n$.
\end{enumerate}
\end{lemma}

\begin{proof}
(i) If $v\le c_n+c'_m$ for all $n,m$, then $v-c_n$ is a lower bound of $(c'_m)$, so $v\le c_n$ for all $n$, and $v\le 0$. A lower bound of $(c_k+c'_k)_k$ is a lower bound of all $c_n+c'_m$, because $c_k+c'_k\le c_n+c'_m$ for $k\ge\max\{n,m\}$.

(ii) If $\pm(z_n-z)\le c_n$ and $\pm(z_n-z')\le c'_n$, then $\pm(z-z')\le c_n+c'_n$ for all $n$, so $\pm(z-z')\le 0$ by (i). Linearity is similar.

(iii) With dominating sequences $(c_n)$ for $(z_n)$ and $(c'_n)$ for $(w_n)$, $z-w=(z-z_n)+(z_n-w_n)+(w_n-w)\le c_n+c'_n$ for large $n$, so $z-w\le 0$ by (i).

(iv) For $m\ge n$, $\pm(z_n-z)\le c_n+c'_m$, where $\pm(z_m-z)\le c'_m\downarrow 0$. So $\pm(z_n-z)-c_n$ is a lower bound of $(c'_m)_{m\ge n}$ and hence $\le 0$.
\end{proof}

We recall the vector metric spaces of \cite{CA09}, allowing the value space to be an arbitrary ordered vector space.

\begin{definition}\label{def:vm}
An \emph{$E$-valued vector metric} on a set $X$ is a map $d\colon X\times X\to E$ such that, for all $x,y,z\in X$,
\begin{enumerate}[label=\textup{(VM\arabic*)}]
\item $d(x,y)=0$ if and only if $x=y$;
\item $d(x,y)\le d(x,z)+d(y,z)$.
\end{enumerate}
A sequence $(x_n)$ \emph{converges} to $x$, written $x_n\to x$, if $d(x_n,x)\le a_n$ for some $a_n\downarrow 0$; it is \emph{Cauchy} if $d(x_n,x_{n+p})\le a_n$ for all $n,p$ and some $a_n\downarrow 0$; and $X$ is \emph{complete} if every Cauchy sequence converges \cite{CA09}.
\end{definition}

The axioms have the expected consequences, with order convergence in $E$ in place of convergence in $\R$.

\begin{lemma}\label{lem:vm}
Let $d$ be an $E$-valued vector metric on $X$.
\begin{enumerate}[label=\textup{(\roman*)}]
\item $d\ge 0$, $d(x,y)=d(y,x)$, and $d(x,y)\le d(x,z)+d(z,y)$.
\item $\pm\bigl(d(x,y)-d(x',y')\bigr)\le d(x,x')+d(y,y')$ for all $x,y,x',y'$.
\item Limits are unique; convergent sequences are Cauchy; and if $x_n\to x$ and $y_n\to y$, then $d(x_n,y_n)\ocv d(x,y)$.
\end{enumerate}
\end{lemma}

\begin{proof}
(i) Taking $y=x$ in (VM2) gives $0\le 2d(x,z)$; taking $z=x$ gives $d(x,y)\le d(y,x)$, and by symmetry of the roles equality holds. The triangle inequality is (VM2) combined with symmetry.

(ii) $d(x,y)\le d(x,x')+d(x',y')+d(y',y)$, and symmetrically.

(iii) If $d(x_n,x)\le a_n$ and $d(x_n,y)\le b_n$, then $d(x,y)\le a_n+b_n$ for all $n$, so $d(x,y)=0$ by Lemma~\ref{lem:olim}(i). Next, $d(x_n,x_{n+p})\le a_n+a_{n+p}\le 2a_n$. Finally, if $d(x_n,x)\le a_n$ and $d(y_n,y)\le b_n$, then $\pm(d(x_n,y_n)-d(x,y))\le a_n+b_n$ by (ii).
\end{proof}

Completions are defined as for metric spaces, with density understood in the sequential sense.

\begin{definition}\label{def:completion}
A \emph{completion} of an $E$-valued vector metric space $(X,d)$ is a complete $E$-valued vector metric space $(\hat X,\hat d)$ together with a map $\iota\colon X\to\hat X$ such that $\hat d(\iota x,\iota y)=d(x,y)$ for all $x,y\in X$ and every point of $\hat X$ is the limit of a sequence in $\iota(X)$.
\end{definition}

We say that convergence in $(X,d)$ satisfies the \emph{iterated limit axiom} if $x^k_n\to x^k$ $(n\to\infty)$ for every $k$ and $x^k\to x$ imply the existence of $n_k$ with $x^k_{n_k}\to x$. The same terminology is used for order convergence of sequences in $E$. Sequential convergence in a first countable topological space satisfies the iterated limit axiom.

The diagonal property of products is governed by a cardinal invariant, which we now recall. For $f,g\in\N^{\N}$ write $f\le^*g$ if $f(k)\le g(k)$ for all but finitely many $k$. A family $\mathcal F\subseteq\N^{\N}$ is \emph{bounded} if some $g$ satisfies $f\le^*g$ for all $f\in\mathcal F$. The \emph{bounding number} $\fb$ is the least cardinality of an unbounded family. It is known that $\aleph_1\le\fb\le 2^{\aleph_0}$, that $\fb=\aleph_1$ under CH, that $\fb=2^{\aleph_0}$ under Martin's axiom, and that $\aleph_1<\fb$ is consistent with ZFC \cite{vD84,Bla10}.

\section{The sequential diagonal property}\label{sec:SDP}

The following two properties of the value space are exactly what completions require (Theorem~\ref{thm:main}).

\begin{definition}\label{def:SOC-SDP}
We say that $E$ has the property
\begin{enumerate}
\item[\SOC] if every order Cauchy sequence in $E$ order converges;
\item[\SDP] if for every family $(a^k_n)$ in $E$ with $a^k_n\downarrow_n0$ for each $k$ there are $n_k\in\N$ and $e_k\downarrow 0$ in $E$ with $a^k_{n_k}\le e_k$ for all $k$. Such a pair $(n_k),(e_k)$ is called a \emph{diagonal} for $(a^k_n)$.
\end{enumerate}
\end{definition}

Diagonals can be thinned out, and \SDP\ need only be tested on families of a special form; this will simplify the constructions in Section~\ref{sec:completions}.

\begin{lemma}\label{lem:basic}
\begin{enumerate}[label=\textup{(\roman*)}]
\item If $(n_k),(e_k)$ is a diagonal for $(a^k_n)$ and $n'_k\ge n_k$ for all $k$, then $(n'_k),(e_k)$ is a diagonal as well.
\item $E$ has \SDP\ if and only if every family with $a^k_n\downarrow_n0$ and $a^k_n\le a^{k+1}_n$ for all $k,n$ has a diagonal.
\item If $E$ is Archimedean and $u\in E_+\setminus\{0\}$, then $E$ has \SDP\ if and only if every family with $a^k_n\downarrow_n0$, $a^k_n\le a^{k+1}_n$ and $a^k_n\ge u/n$ has a diagonal.
\end{enumerate}
\end{lemma}

\begin{proof}
(i) $a^k_{n'_k}\le a^k_{n_k}$. (ii) Given $(a^k_n)$, the family $b^k_n=\sum_{i\le k}a^i_n$ satisfies $b^k_n\downarrow_n0$ by Lemma~\ref{lem:olim}(i), increases in $k$ and dominates $a^k_n$. (iii) Replace $b^k_n$ by $b^k_n+u/n$, which still decreases to $0$ because $E$ is Archimedean.
\end{proof}

In terms of convergence, \SDP\ is the iterated limit axiom.

\begin{proposition}\label{prop:iterated}
The space $E$ has \SDP\ if and only if order convergence of sequences in $E$ satisfies the iterated limit axiom.
\end{proposition}

\begin{proof}
Assume \SDP, $\pm(z^k_n-z^k)\le a^k_n\downarrow_n0$ and $\pm(z^k-z)\le b_k\downarrow 0$. For a diagonal $(n_k),(e_k)$ of $(a^k_n)$ we get $\pm(z^k_{n_k}-z)\le e_k+b_k\downarrow 0$. Conversely, given $(a^k_n)$, apply the iterated limit axiom to $z^k_n=a^k_n$, $z^k=z=0$.
\end{proof}

A first consequence concerns positive functionals. It gives a convenient test for the failure of \SDP; its converse is discussed in Section~\ref{sec:further}.

\begin{proposition}\label{prop:functionals}
If $E$ has \SDP, every positive linear functional on $E$ is $\sigma$-order continuous, i.e.\ $x_n\downarrow 0$ implies $f(x_n)\to 0$.
\end{proposition}

\begin{proof}
Let $f\ge 0$, $x_n\downarrow 0$, and suppose $f(x_n)\downarrow c>0$. A diagonal for $a^k_n=kx_n$ gives $kx_{n_k}\le e_k\le e_1$, so $kc\le f(e_1)$ for all $k$, which is impossible.
\end{proof}

Order Cauchy completeness, in contrast, is automatic when countable suprema exist.

\begin{proposition}\label{prop:dedekind}
Every Dedekind $\sigma$-complete vector lattice, that is, one in which every nonempty countable set that is bounded above has a supremum, satisfies \SOC.
\end{proposition}

\begin{proof}
Let $\pm(z_n-z_m)\le c_n$ for $m\ge n$, where $c_n\downarrow 0$. The sequence lies in $[z_1-c_1,z_1+c_1]$, so $z=\inf_n\sup_{m\ge n}z_m$ exists. For $m\ge n$, $z_m\le z_n+c_n$, hence $z\le z_n+c_n$; and $\sup_{m\ge N}z_m\ge z_N\ge z_n-c_n$ for $N\ge n$, hence $z\ge z_n-c_n$.
\end{proof}

We now compare \SDP\ with the classical properties of Egoroff type. The following conditions are classical for vector lattices; see Remark~\ref{rem:classical}.

\begin{definition}\label{def:egoroff}
An ordered vector space $E$ has
\begin{enumerate}[label=\textup{(\roman*)}]
\item the \emph{dominating property} if for every family $(a^k_n)$ with $a^k_n\downarrow_n0$ for each $k$ there are $n_k\in\N$ and $w\in E_+$ with $a^k_{n_k}\le w$ for all $k$;
\item the \emph{$\sigma$-property} if for every sequence $(u_n)$ in $E_+$ there are reals $\lambda_n>0$ and $w\in E_+$ with $\lambda_nu_n\le w$ for all $n$;
\item property \OCRU\ if every sequence $c_n\downarrow 0$ converges to $0$ \emph{relatively uniformly}, i.e.\ there are $w\in E_+$ and reals $\eps_n\ge 0$ with $\eps_n\to 0$ and $c_n\le\eps_nw$ for all $n$.
\end{enumerate}
The space $E$ is \emph{uniformly complete} if for every $w\in E_+$ and every sequence $(z_n)$ with $\pm(z_n-z_m)\le\eps_nw$ for $m\ge n$, where $\eps_n\ge 0$ and $\eps_n\to 0$, there are $z\in E$, $v\in E_+$ and reals $\delta_n\to 0$ with $\pm(z_n-z)\le\delta_nv$ for all $n$.
\end{definition}

Their relation to \SDP\ is the same as in vector lattices, but the proofs have to avoid lattice operations.

\begin{proposition}\label{prop:egoroff}
For an Archimedean ordered vector space $E$ the following are equivalent:
\begin{enumerate}[label=\textup{(\alph*)}]
\item $E$ has \SDP;
\item $E$ has the dominating property;
\item $E$ has the $\sigma$-property and property \OCRU.
\end{enumerate}
\end{proposition}

\begin{proof}
(a)$\Rightarrow$(b). If $(n_k),(e_k)$ is a diagonal, then $a^k_{n_k}\le e_k\le e_1$.

(b)$\Rightarrow$(a). Given $(a^k_n)$, apply (b) to the family $(ka^k_n)$, which also satisfies $ka^k_n\downarrow_n0$. This gives $w$ and $n_k$ with $a^k_{n_k}\le w/k$, and $w/k\downarrow 0$ because $E$ is Archimedean.

(a), (b)$\Rightarrow$(c). Let $(u_n)$ be a sequence in $E_+$. For each $n$, $u_n/k\downarrow_k0$, so (b) gives $k_n$ and $w$ with $u_n/k_n\le w$; take $\lambda_n=1/k_n$. Now let $c_n\downarrow 0$. A diagonal $(n_k),(e_k)$ of the family $(kc_n)$ gives $kc_{n_k}\le e_1$, and by Lemma~\ref{lem:basic}(i) we may assume $n_1<n_2<\cdots$. With $w=c_1+e_1$, $\eps_n=1$ for $n<n_1$ and $\eps_n=1/k$ for $n_k\le n<n_{k+1}$, we get $c_n\le c_{n_k}\le\eps_nw$ for $n\ge n_1$, and $c_n\le c_1\le w$ for $n<n_1$.

(c)$\Rightarrow$(a). Given $(a^k_n)$, \OCRU\ gives $w_k\in E_+$ and $\eps^k_n\to 0$ $(n\to\infty)$ with $a^k_n\le\eps^k_nw_k$, and the $\sigma$-property gives $\lambda_k>0$ and $w\in E_+$ with $\lambda_kw_k\le w$. Choose $n_k$ with $\eps^k_{n_k}\le\lambda_k/k$. Then $a^k_{n_k}\le w/k$, and $w/k\downarrow 0$.
\end{proof}

For uniformly complete spaces, \SOC\ is a consequence of \SDP; we use this for spaces of continuous functions in Corollary~\ref{cor:CX}.

\begin{lemma}\label{lem:SOCuniform}
If $E$ is Archimedean, uniformly complete and has \SDP, then $E$ satisfies \SOC.
\end{lemma}

\begin{proof}
Let $\pm(z_n-z_m)\le c_n$ for $m\ge n$ with $c_n\downarrow 0$. By Proposition~\ref{prop:egoroff}, $c_n\le\eps_nw$ with $\eps_n\to 0$, so there are $z$, $v$ and $\delta_n\to 0$ with $\pm(z_n-z)\le\delta_nv$. With $\delta'_n=\sup_{m\ge n}|\delta_m|$ we have $\pm(z_n-z)\le\delta'_nv$, and $\delta'_nv\downarrow 0$ because $E$ is Archimedean.
\end{proof}

\begin{remark}[the classical case]\label{rem:classical}
Let $E$ be a vector lattice.
\begin{enumerate}[label=(\alph*)]
\item Property \SDP\ is the \emph{diagonal property for order convergence}, which goes back to Kantorovich \cite[Satz~24]{Kan37}; by Proposition~\ref{prop:iterated} it says that order convergence of sequences satisfies the iterated limit axiom. Dodds \cite[Theorem~3.9]{Dod74} proved that it is equivalent to the \emph{strong Egoroff property}: whenever $u_{n,k}\downarrow_k0$ for every $n$, there are $v_m\downarrow 0$ such that for all $m,n$ some $k$ satisfies $u_{n,k}\le v_m$. It is also equivalent to the idempotence $(S')'=S'$ of the pseudo order closure $S'$, the set of order limits of sequences in $S$, for every $S\subseteq E$ \cite[Theorems~3.9 and 3.11]{Dod74}; see also \cite[Theorems~68.8 and 68.9]{LZ71}. Our single-index formulation is equivalent to the double-index one in every ordered vector space: given $(a^k_n)$, take a diagonal $(n_k),(e_k)$ of the family $b^k_n=\sum_{i\le k}a^i_n$ from the proof of Lemma~\ref{lem:basic}(ii); for all $m,k$ and $l=\max\{m,k\}$ we get $e_m\ge e_l\ge b^l_{n_l}\ge a^k_{n_l}$. Conversely, applying the strong Egoroff property to $u_{n,k}=a^n_k$ and taking $m=n$ gives indices $k_n$ with $a^n_{k_n}\le v_n$, that is, a diagonal.
\item The \emph{Egoroff property} is the corresponding condition for order bounded families. It is strictly weaker; for instance, $\ell^\infty$ has it but fails \SDP\ \cite[p.~368]{Dod74}. The strong Egoroff property is the conjunction of the Egoroff property and the dominating property \cite[Theorem~3.4]{Dod74}.
\item For Archimedean vector lattices, (a)$\Leftrightarrow$(b) in Proposition~\ref{prop:egoroff} is \cite[Theorem~3.5]{Dod74} and \cite[Theorem~68.5]{LZ71}; the dominating property is the d-property of \cite[Definition~68.3]{LZ71}, which goes back to Nakano; see also \cite[p.~372]{Dod74}. The equivalence (a)$\Leftrightarrow$(c) is \cite[Theorem~70.2]{LZ71}, where \OCRU\ appears as stability of order convergence; see also \cite{Zaa68}. Condition (c) is the definition of a \emph{strongly Egoroff} vector lattice in \cite{HvM15}. The $\sigma$-property alone characterises the diagonal property for relatively uniform convergence \cite[Theorem~4.6]{Dod74}, \cite[Theorem~72.2]{LZ71}.
\end{enumerate}
The proofs of Propositions~\ref{prop:iterated} and \ref{prop:egoroff} use no lattice operations; this is what we need for the ordered vector spaces and pre-Riesz spaces in Sections~\ref{sec:banach} and \ref{sec:perm}.
\end{remark}

\section{Completions}\label{sec:completions}

We begin with the sufficiency of the two properties. The construction is the classical one by classes of Cauchy sequences: \SOC\ is needed to define the distance on the quotient, and \SDP\ to prove that the quotient is complete.

\begin{theorem}\label{thm:exist}
If $E$ satisfies \SOC\ and \SDP, then every $E$-valued vector metric space $(X,d)$ has a completion.
\end{theorem}

\begin{proof}
\emph{Construction.} Let $\mathcal C$ be the set of Cauchy sequences in $X$. For $\mathbf x=(x_n)$ and $\mathbf y=(y_n)$ in $\mathcal C$ with Cauchy bounds $(a_n)$ and $(b_n)$, that is, $d(x_n,x_{n+p})\le a_n$ and $d(y_n,y_{n+p})\le b_n$ for all $n,p$ with $a_n,b_n\downarrow 0$, Lemma~\ref{lem:vm}(ii) gives $\pm(d(x_n,y_n)-d(x_m,y_m))\le a_n+b_n$ for $m\ge n$. So $(d(x_n,y_n))$ is order Cauchy, and by \SOC\ it has an order limit $D(\mathbf x,\mathbf y)$. By Lemma~\ref{lem:olim}, $D$ is symmetric, satisfies (VM2), and $D(\mathbf x,\mathbf x)=0$. Hence $\mathbf x\sim\mathbf y\Leftrightarrow D(\mathbf x,\mathbf y)=0$ is an equivalence relation. By Lemma~\ref{lem:vm}(ii), $D(\mathbf x,\mathbf y)=D(\mathbf x',\mathbf y')$ whenever $\mathbf x\sim\mathbf x'$ and $\mathbf y\sim\mathbf y'$, so $\hat d([\mathbf x],[\mathbf y])=D(\mathbf x,\mathbf y)$ is a well-defined vector metric on $\hat X=\mathcal C/{\sim}$. For $x\in X$ let $\iota(x)$ be the class of the constant sequence $(x,x,\dots)$; then $\hat d(\iota x,\iota y)=d(x,y)$.

\emph{Density.} Let $\xi=[(x_n)]$ with Cauchy bound $(a_n)$. Since $d(x_m,x_n)\le a_n$ for $m\ge n$, Lemma~\ref{lem:olim}(iii) gives $\hat d(\xi,\iota x_n)\le a_n$, so $\iota(x_n)\to\xi$.

\emph{Completeness.} Let $(\xi_k)$ be Cauchy in $\hat X$, say $\hat d(\xi_k,\xi_{k+p})\le b_k$ for all $k,p$ with $b_k\downarrow 0$. Choose representatives $(x^k_n)_n$ of $\xi_k$ with Cauchy bounds $a^k_n\downarrow_n0$, so that $\hat d(\xi_k,\iota x^k_n)\le a^k_n$. By \SDP\ there is a diagonal $(n_k),(e_k)$. With $y_k=x^k_{n_k}$,
\[
d(y_k,y_{k+p})\le\hat d(\iota y_k,\xi_k)+\hat d(\xi_k,\xi_{k+p})+\hat d(\xi_{k+p},\iota y_{k+p})\le 2e_k+b_k,
\]
so $(y_k)\in\mathcal C$ with Cauchy bound $2e_k+b_k$. For $\xi=[(y_k)]$ the density step gives $\hat d(\iota y_k,\xi)\le 2e_k+b_k$, hence $\hat d(\xi_k,\xi)\le 3e_k+b_k\downarrow 0$.
\end{proof}

As for metric spaces, the completion is determined by the space.

\begin{theorem}\label{thm:unique}
Completions are unique up to isometry: if $(\hat X_1,\hat d_1,\iota_1)$ and $(\hat X_2,\hat d_2,\iota_2)$ are completions of $(X,d)$, then there is a unique bijection $\Phi\colon\hat X_1\to\hat X_2$ with $\Phi\circ\iota_1=\iota_2$ and $\hat d_2(\Phi\xi,\Phi\eta)=\hat d_1(\xi,\eta)$.
\end{theorem}

\begin{proof}
For $\xi\in\hat X_1$ choose a sequence $(x_n)$ in $X$ with $\iota_1(x_n)\to\xi$. Then $(x_n)$ is Cauchy in $X$, so $(\iota_2x_n)$ converges in $\hat X_2$ to some $\eta$; put $\Phi(\xi)=\eta$. If $\iota_1(x_n)\to\xi$ and $\iota_1(x'_n)\to\xi'$, then by Lemma~\ref{lem:vm}(iii)
\[
\hat d_1(\xi,\xi')\xleftarrow{\ \mathrm o\ }d(x_n,x'_n)\xrightarrow{\ \mathrm o\ }\hat d_2(\Phi\xi,\Phi\xi'),
\]
so the two limits agree. For $\xi=\xi'$ this shows that $\Phi$ is well defined, and in general that $\Phi$ is an isometry. Exchanging the roles gives the inverse. Uniqueness holds because an isometry preserves convergence and $\iota_1(X)$ is dense.
\end{proof}

The necessity part of Theorem~\ref{thmA} uses a family of countable vector metric spaces.

\begin{lemma}\label{lem:rows}
Let $E$ be an ordered vector space, $(b^k_n)$ a family in $E_+\setminus\{0\}$ with $b^k_n$ decreasing in $n$, and $(w_k)$ a strictly decreasing sequence in $E_+\setminus\{0\}$; put $w_\infty=0$. On the set
\[
\bar X=\{p^k_n:k,n\in\N\}\cup\{\xi_k:k\in\N\}\cup\{\xi\}
\]
define a \emph{row} $r$ and a \emph{height} $h$ by $r(p^k_n)=r(\xi_k)=k$, $r(\xi)=\infty$, $h(p^k_n)=b^k_n$ and $h(\xi_k)=h(\xi)=0$. For distinct $P,Q$ let
\[
d(P,Q)=\begin{cases}\max\{h(P),h(Q)\}, & r(P)=r(Q),\\ h(P)+h(Q)+w_{\min\{r(P),r(Q)\}}, & r(P)\ne r(Q),\end{cases}
\]
and $d(P,P)=0$. Then $d$ is an $E$-valued vector metric on $\bar X$, and its restriction to $X_0=\{p^k_n:k,n\in\N\}$ is a vector metric on $X_0$. Moreover, $p^k_n\to\xi_k$ as $n\to\infty$ whenever $b^k_n\downarrow_n0$, and $\xi_k\to\xi$ whenever $w_k\downarrow 0$.
\end{lemma}

\begin{proof}
Within a row the heights form a chain ($0\le\dots\le b^k_2\le b^k_1$), so any two heights are comparable and their maximum exists.

(VM1). Two distinct points in a finite row cannot both have height $0$, and row $\infty$ has only one point; two points in different rows have distance at least $w_{\min\{r(P),r(Q)\}}>0$, as one of the rows is finite.

(VM2). We verify $d(P,Q)\le d(P,R)+d(Q,R)$ for distinct $P,Q,R$ by cases, writing $h_P=h(P)$.
\begin{enumerate}[label=(\alph*)]
\item All three in one row: $\max\{h_P,h_Q\}\le\max\{h_P,h_R\}+\max\{h_Q,h_R\}$, since all heights in a row are comparable and positive.
\item $P,Q$ in one row, $R$ in another: $\max\{h_P,h_Q\}\le h_P+h_Q$, which is at most the right-hand side.
\item $P,R$ in one row, $Q$ in another: the claim reads $h_P+h_Q+w\le\max\{h_P,h_R\}+h_R+h_Q+w$ with the same $w$, i.e.\ $h_P\le\max\{h_P,h_R\}+h_R$. The case of $Q,R$ in one row is symmetric.
\item Three different rows $k,l,j$: the claim reduces to $w_{\min\{k,l\}}\le w_{\min\{k,j\}}+w_{\min\{l,j\}}$. Since $w$ is nonincreasing on $\N\cup\{\infty\}$ and $\min\{k,j\}\le\min\{k,l\}$ or $\min\{l,j\}\le\min\{k,l\}$, one of the two terms on the right already dominates the left.
\end{enumerate}
Finally, $d(p^k_n,\xi_k)=b^k_n$ and $d(\xi_k,\xi)=w_k$.
\end{proof}

For the rest of this section, $E\ne\{0\}$ is Archimedean and $u\in E_+\setminus\{0\}$; note that $u/n\downarrow 0$.

The row spaces first show that \SDP\ is the iterated limit axiom for all $E$-valued vector metric spaces.

\begin{proposition}\label{prop:iterated-vm}
If convergence in every $E$-valued vector metric space satisfies the iterated limit axiom, then $E$ has \SDP. Conversely, \SDP\ implies the iterated limit axiom in every $E$-valued vector metric space.
\end{proposition}

\begin{proof}
Let $(a^k_n)$ be given; by Lemma~\ref{lem:basic}(iii) we may assume $a^k_n\ge u/n$ and $a^k_n\le a^{k+1}_n$. Apply Lemma~\ref{lem:rows} with $b^k_n=a^k_n$ and $w_k=u/k$. Then $p^k_n\to\xi_k$ and $\xi_k\to\xi$, so by assumption there are $n_k$ and $e_k\downarrow 0$ with $d(p^k_{n_k},\xi)=a^k_{n_k}+w_k\le e_k$. Hence $a^k_{n_k}\le e_k$.

Conversely, let $d(x^k_n,x^k)\le a^k_n\downarrow_n0$ and $d(x^k,x)\le b_k\downarrow 0$. A diagonal $(n_k),(e_k)$ for $(a^k_n)$ gives $d(x^k_{n_k},x)\le e_k+b_k\downarrow 0$.
\end{proof}

In particular, when \SDP\ fails, convergence in vector metric spaces cannot always be described by a metric.

\begin{corollary}\label{cor:nontop}
If $E$ fails \SDP, then there is a countable $E$-valued vector metric space whose convergence is not induced by any first countable topology, in particular by no metric.
\end{corollary}

\begin{proof}
If \SDP\ fails, the proof of Proposition~\ref{prop:iterated-vm} produces a family for which the countable space $\bar X$ of Lemma~\ref{lem:rows} violates the iterated limit axiom, whereas convergence of sequences in a first countable topology satisfies it.
\end{proof}

We now prove that both properties are necessary for the existence of completions, beginning with \SOC.

\begin{proposition}\label{prop:SOCnec}
If every $E$-valued vector metric space has a completion, then $E$ satisfies \SOC.
\end{proposition}

\begin{proof}
Let $\pm(z_n-z_m)\le c_n$ for $m\ge n$, where $c_n\downarrow 0$. Put $z'_n=z_n-z_1+c_1\ge 0$, $a_n=c_n+u/n$ (so $a_n\downarrow 0$ and $a_n>0$) and $v=a_1$. On $X=\{q\}\cup\{p_n:n\in\N\}$ let $d(p_n,p_m)=a_{\min\{n,m\}}$ for $n\ne m$, $d(p_n,q)=d(q,p_n)=z'_n+v$, and $d=0$ on the diagonal. (VM1) is clear. For (VM2): among the $p$'s use that $(a_n)$ decreases; $d(p_n,p_m)\le v\le d(p_n,q)+d(p_m,q)$; and $d(p_n,q)\le d(p_n,p_m)+d(p_m,q)$ because $z'_n-z'_m=z_n-z_m\le c_{\min\{n,m\}}\le a_{\min\{n,m\}}$.

Since $d(p_n,p_{n+r})=a_n$, $(p_n)$ is Cauchy. In a completion, $\iota(p_n)\to\zeta$ for some $\zeta$, and $z'_n+v=\hat d(\iota p_n,\iota q)\ocv\hat d(\zeta,\iota q)$ by Lemma~\ref{lem:vm}(iii). Hence $(z_n)$ order converges.
\end{proof}

The necessity of \SDP\ is the main step: a completion of a suitable row space produces a diagonal.

\begin{proposition}\label{prop:SDPnec}
If every $E$-valued vector metric space has a completion, then $E$ satisfies \SDP.
\end{proposition}

\begin{proof}
Let $(a^k_n)$ be given; as before we may assume $a^k_n\ge u/n$ and $a^k_n\le a^{k+1}_n$. Let $b^k_n=a^k_n$, $w_k=u/k$, and let $X_0=\{p^k_n\}$ carry the vector metric of Lemma~\ref{lem:rows}. Let $(\hat X,\hat d,\iota)$ be a completion; we identify $p^k_n$ with $\iota(p^k_n)$.

Each row $(p^k_n)_n$ is Cauchy, since $d(p^k_n,p^k_{n+r})=b^k_n$, so $p^k_n\to\xi_k$ for some $\xi_k\in\hat X$. Lemma~\ref{lem:vm}(iii) and uniqueness of order limits give, for $k\ne l$,
\[
\hat d(\xi_k,p^k_n)=b^k_n,\qquad\hat d(\xi_l,p^k_n)=b^k_n+w_{\min\{k,l\}},\qquad\hat d(\xi_k,\xi_l)=w_{\min\{k,l\}}.
\]
For instance, $\hat d(\xi_k,p^k_n)$ is the order limit, as $m\to\infty$, of $d(p^k_m,p^k_n)=b^k_n$ $(m>n)$; and $\hat d(\xi_k,\xi_l)$ is the order limit of $d(p^k_m,p^l_m)=b^k_m+b^l_m+w_{\min\{k,l\}}$. So $(\xi_k)$ is Cauchy and converges to some $\xi$, with $\hat d(\xi_k,\xi)=\lim_lw_{\min\{k,l\}}=w_k$. By density there are $y_j=p^{k_j}_{n_j}$ and $e'_j\downarrow 0$ with $\hat d(\xi,y_j)\le e'_j$.

\emph{Each $k$ occurs only finitely often among the $k_j$.} Suppose $k_j=k$ for all $j$ in an infinite set $J$. If some $n$ satisfies $n_j=n$ for infinitely many $j\in J$, then $\hat d(\xi,p^k_n)\le e'_j$ for those $j$, so $\xi=p^k_n$, and for $l>k$ we get $w_l=\hat d(\xi_l,\xi)=\hat d(\xi_l,p^k_n)=b^k_n+w_k>w_l$, a contradiction. Otherwise, for all $N,M$ there is $j\in J$ with $n_j\ge N$ and $j\ge M$, and then $w_k=\hat d(\xi_k,\xi)\le\hat d(\xi_k,y_j)+\hat d(y_j,\xi)\le b^k_N+e'_M$. By Lemma~\ref{lem:olim}(i), $w_k\le 0$, a contradiction.

\emph{The diagonal.} Let $J(k)=1+\max\{j:k_j<k\}$, with $\max\emptyset=0$. Then $k_j\ge k$ for $j\ge J(k)$, $J(k)$ is nondecreasing, and $J(k)\to\infty$. Put $m_k=n_{J(k)}$, $k'=k_{J(k)}\ge k$, and $c_k=w_k+e'_{J(k)}$; then $c_k\downarrow 0$ by Lemma~\ref{lem:olim}(i). Since $b^k_n\le b^{k'}_n$,
\[
a^k_{m_k}=b^k_{m_k}\le b^{k'}_{m_k}=\hat d(\xi_{k'},y_{J(k)})\le\hat d(\xi_{k'},\xi)+\hat d(\xi,y_{J(k)})\le w_{k'}+e'_{J(k)}\le c_k.
\]
Thus $(m_k),(c_k)$ is a diagonal for $(a^k_n)$.
\end{proof}

Collecting the results of this section, we obtain the main theorem, which contains Theorem~\ref{thmA}.

\begin{theorem}\label{thm:main}
Let $E\ne\{0\}$ be an Archimedean ordered vector space. The following are equivalent:
\begin{enumerate}[label=\textup{(\alph*)}]
\item every $E$-valued vector metric space has a completion;
\item every countable $E$-valued vector metric space has a completion;
\item $E$ satisfies \SOC\ and \SDP.
\end{enumerate}
The completion is unique up to isometry. Moreover, $E$ has \SDP\ if and only if convergence in every (countable) $E$-valued vector metric space satisfies the iterated limit axiom.
\end{theorem}

\begin{proof}
(a)$\Rightarrow$(b) is trivial; (b)$\Rightarrow$(c) by Propositions~\ref{prop:SOCnec} and \ref{prop:SDPnec}, whose spaces are countable; (c)$\Rightarrow$(a) by Theorem~\ref{thm:exist}. Uniqueness is Theorem~\ref{thm:unique}, and the last statement is Proposition~\ref{prop:iterated-vm}.
\end{proof}

\begin{remark}\label{rem:closure}
For a vector lattice $E$, the last statement of Theorem~\ref{thm:main} and Corollary~\ref{cor:nontop} are vector metric counterparts of Dodds' theorem that \SDP\ is equivalent to the idempotence of the pseudo order closure in $E$ \cite[Theorems~3.9 and 3.11]{Dod74}. The spaces of Lemma~\ref{lem:rows} make a failure of \SDP\ visible inside a single countable vector metric space.
\end{remark}

\begin{remark}\label{rem:Cev15}
Completions of vector metric spaces were constructed in \cite[Theorem~4]{Cev15} for order Cauchy complete vector lattices $E$, that is, under \SOC. There, completeness of $\hat X$ is obtained by choosing, for a Cauchy sequence $(\hat x_n)$ in $\hat X$, points $y_n\in X$ with $\hat d(\iota(y_n),\hat x_n)\le a_n$ for one sequence $a_n\downarrow 0$. Density provides such estimates for each $n$ separately, and combining them is a diagonal choice; this is exactly where \SDP\ enters, as in the proof of Theorem~\ref{thm:exist}. Example~\ref{ex:linfty} shows that the step cannot be carried out under \SOC\ alone, and Theorem~\ref{thm:main} identifies \SDP\ as the missing condition: the construction of \cite{Cev15} works for every $E$-valued vector metric space precisely when $E$ satisfies \SOC\ and \SDP. The same diagonal step occurs in the construction of unbounded completions in \cite[Theorem~2]{OC20}. For $E=\R^{\fb}$, unbounded order convergence of sequences coincides with order convergence, so Corollary~\ref{cor:RA} and Proposition~\ref{prop:SDPnec} show that not every $\R^{\fb}$-valued vector metric space has an unbounded completion, although $\R^{\fb}$ is uo-Cauchy complete; an unbounded analogue of \SDP\ is needed there.
\end{remark}
\section{Ordered Banach spaces and the spaces \texorpdfstring{$C(K)$}{C(K)}}\label{sec:banach}

An ordered Banach space is a real Banach space with a closed cone $E_+$. The cone is \emph{normal} if there is a constant $C$ with $\|x\|\le C\|y\|$ whenever $0\le x\le y$; then $-c\le x\le c$ implies $\|x\|\le(2C+1)\|c\|$. The norm is \emph{$\sigma$-order continuous} if $z_n\downarrow 0$ implies $\|z_n\|\to 0$.

For these spaces the diagonal property turns out to be a property of the norm.

\begin{theorem}\label{thm:banach}
Let $E$ be an ordered Banach space with a closed normal cone. The following are equivalent:
\begin{enumerate}[label=\textup{(\alph*)}]
\item $E$ has \SDP;
\item every continuous positive linear functional on $E$ is $\sigma$-order continuous;
\item the norm of $E$ is $\sigma$-order continuous.
\end{enumerate}
If they hold, $E$ also satisfies \SOC, and hence every $E$-valued vector metric space has a completion.
\end{theorem}

\begin{proof}
(a)$\Rightarrow$(b) is Proposition~\ref{prop:functionals}.

(b)$\Rightarrow$(c). As $E_+$ is closed and normal, the dual cone generates the dual space (Krein--Grosberg duality, see \cite[Theorem~2.40 and Corollary~2.43]{AT07}). Hence (b) implies that $z_n\downarrow 0$ entails $f(z_n)\to 0$ for every continuous linear $f$, i.e.\ $z_n\to 0$ weakly. By Mazur's theorem, for each $m$ there is a convex combination $y_m$ of $z_m,\dots,z_{N_m}$ with $\|y_m\|<1/m$. Since $(z_n)$ decreases, $0\le z_j\le y_m$ for $j\ge N_m$, and normality gives $\|z_j\|\le C/m$.

(c)$\Rightarrow$(a). Given $(a^k_n)$, choose $n_k$ with $\|a^k_{n_k}\|\le 2^{-k}$ and let $e_k=\sum_{j\ge k}a^j_{n_j}$. Since $E_+$ is closed, $e_k\ge a^k_{n_k}$ and $(e_k)$ decreases; as $\|e_k\|\to 0$, every lower bound $h$ satisfies $h\le 0$, so $e_k\downarrow 0$.

\SOC: if $\pm(z_n-z_m)\le c_n$ for $m\ge n$, then $\|z_n-z_m\|\le(2C+1)\|c_n\|\to 0$, so $z_n\to z$ in norm, and $\pm(z_n-z)\le c_n$ because $E_+$ is closed. The last assertion is Theorem~\ref{thm:exist}.
\end{proof}

Normality is used only for (b)$\Rightarrow$(c) and for \SOC; (c)$\Rightarrow$(a) holds for every ordered Banach space with a closed cone.

\begin{corollary}\label{cor:BL}
A Banach lattice has \SDP\ if and only if its norm is $\sigma$-order continuous. In that case every vector metric space with values in it has a completion.
\end{corollary}

For Banach lattices, Corollary~\ref{cor:BL} can also be derived from the classical theory. A Banach lattice has the $\sigma$-property, with $w=\sum_n2^{-n}(1+\|u_n\|)^{-1}u_n$. Its norm is $\sigma$-order continuous if and only if it has \OCRU: if $u_n\downarrow 0$ and $\|u_n\|\to 0$, choose $n_1<n_2<\cdots$ with $\|u_{n_k}\|\le 4^{-k}$ and put $w=u_1+\sum_k2^ku_{n_k}$; then $u_n\le 2^{-k}w$ for $n\ge n_k$. Conversely, $u_n\le\eps_nw$ gives $\|u_n\|\le\eps_n\|w\|$. So Corollary~\ref{cor:BL} follows from Proposition~\ref{prop:egoroff}, or from \cite[Theorem~70.2]{LZ71}. Theorem~\ref{thm:banach} is needed for spaces that are not lattices, such as those in Examples~\ref{ex:good}(b), \ref{ex:spin} and \ref{ex:vN}.

\begin{example}[admissible value spaces]\label{ex:good}
\begin{enumerate}[label=(\alph*)]
\item $L^p(\mu)$ for $1\le p<\infty$, $c_0(\Gamma)$ for every set $\Gamma$, and every finite-dimensional space with a closed cone (closed cones in finite dimensions are normal). In the last case, $z_n\downarrow 0$ implies $\|z_n\|\to 0$ because $(z_n)$ lies in the compact order interval $[0,z_1]$ and every cluster point is a lower bound.
\item The space $S_2^{\mathrm{sa}}$ of self-adjoint Hilbert--Schmidt operators on a separable Hilbert space with the L\"owner cone. It is not a vector lattice: two incomparable operators $A,B$ supported on a common two-dimensional subspace $V$ have no supremum. Indeed, if $S$ were their supremum and $P$ the projection onto $V$, then $PSP$ would also be an upper bound, so $T=PSP-S\ge 0$; moreover $\langle Tx,x\rangle=0$ for $x\in V$ and for $x\perp V$ (as $\langle Sx,x\rangle\ge\langle Ax,x\rangle=0$ there), so $T=0$ and $S=PSP$. Then $S|_V$ would be a supremum of $A|_V$ and $B|_V$ in the self-adjoint $2\times 2$ matrices, contradicting Kadison's theorem that these form an anti-lattice \cite{Kad51}. The cone is closed and normal, since $0\le A\le B$ implies $\operatorname{tr}A^2\le\operatorname{tr}(AB)\le\operatorname{tr}B^2$. The norm is $\sigma$-order continuous: if $A_n\downarrow 0$, then $A_n$ converges strongly to some $A$ with $0\le A\le A_n$, so $A\in S_2^{\mathrm{sa}}$ is a lower bound and $A=0$; and $\operatorname{tr}A_n^2\le\operatorname{tr}(A_1A_n)=\sum_i\lambda_i\langle A_ne_i,e_i\rangle\to 0$ by dominated convergence, where $(e_i)$ is an orthonormal eigenbasis of $A_1$ with eigenvalues $(\lambda_i)$.
\item $C(K)$ for $K=D\cup\{\infty\}$, the one-point compactification of an uncountable discrete set $D$. Let $x_n\downarrow 0$ in $C(K)$. Since $\mathbf 1_{\{d\}}\in C(K)$, $\inf_nx_n(d)=0$ for each $d\in D$. If $x_n(\infty)\ge c>0$ for all $n$, then $x_n\ge c/2$ off a finite set $F_n\subseteq D$, and any $d\notin\bigcup_nF_n$ gives a contradiction. So $x_n(\infty)\to 0$. If $\|x_n\|_\infty\ge c$ for all $n$, pick $d_n\in D$ with $x_n(d_n)\ge c$ for large $n$. Infinitely many $d_n$ are distinct, and these converge to $\infty$. For fixed $m$, $x_m(d_n)\ge x_n(d_n)\ge c$ for $n\ge m$, so $x_m(\infty)\ge c$, a contradiction. For countable $D$, $C(K)\cong c$ fails \SDP; see Example~\ref{ex:bad}. For uncountable $D$, $K$ is an almost $P$-space \cite[Remark~6.5(c)]{HvM15}, so this is also a case of Theorem~\ref{thm:CK}.
\end{enumerate}
\end{example}

On the other hand, the most familiar Banach lattices with an order unit fail \SDP.

\begin{example}[spaces failing \SDP]\label{ex:bad}
The spaces $\ell^\infty$, $c$ (convergent sequences) and $C[0,1]$ fail \SDP\ by Corollary~\ref{cor:BL}: $\mathbf 1_{[n,\infty)}\downarrow 0$ in $\ell^\infty$ and in $c$, and tents of height $1$ at $\frac12$ with shrinking supports decrease to $0$ in $C[0,1]$, all of norm $1$. Alternatively, Proposition~\ref{prop:functionals} applies: Banach limits, the limit functional on $c$ and the point evaluation at $\frac12$ are positive functionals that are not $\sigma$-order continuous.
\end{example}

For $\ell^\infty$ the failure can be made explicit; the example also shows that Dedekind completeness of the value space does not guarantee completions.

\begin{example}[an $\ell^\infty$-valued space without completion]\label{ex:linfty}
The space $\ell^\infty$ is Dedekind complete, so it satisfies \SOC\ by Proposition~\ref{prop:dedekind}, but it fails \SDP. The family $a^k_n=k\mathbf 1_{[n,\infty)}$ has no diagonal, since $k\mathbf 1_{[n_k,\infty)}\le e_k\le e_1$ for all $k$ would make $e_1$ unbounded. Following the proof of Proposition~\ref{prop:SDPnec} with $u=\mathbf 1$, let
\[
b^k_n=\tfrac{k(k+1)}{2}\mathbf 1_{[n,\infty)}+\tfrac1n\mathbf 1,\qquad w_k=\tfrac1k\mathbf 1,
\]
and $X_0=\{p^k_n:k,n\in\N\}$ with $d(p^k_n,p^k_m)=b^k_{\min\{n,m\}}$ for $n\ne m$ and $d(p^k_n,p^l_m)=b^k_n+b^l_m+w_{\min\{k,l\}}$ for $k\ne l$. This is an explicit countable $\ell^\infty$-valued vector metric space without a completion. Regarded as an $\R^{\N}$-valued vector metric space, the same $X_0$ has a completion, by Corollary~\ref{cor:RA} below. Completions therefore depend on the value space, not only on $X_0$ and the values of $d$.
\end{example}

If order convergence and norm convergence of sequences coincide in $E$, then convergence in every $E$-valued vector metric space is described by the scalar metric $\|d\|$. For cones with interior points this is again a matter of $\sigma$-order continuity. Recall that $u\in E_+$ is an \emph{order unit} if for every $x\in E$ there is $\lambda\ge 0$ with $-\lambda u\le x\le\lambda u$; the \emph{order unit norm} is $\|x\|_u=\inf\{\lambda\ge 0:-\lambda u\le x\le\lambda u\}$, and every interior point of a closed cone is an order unit.

\begin{proposition}\label{prop:unit}
Let $E$ be an ordered Banach space whose closed cone is normal and has an interior point $u$. Then every norm null sequence is order null. Hence order convergence of sequences coincides with norm convergence if and only if the norm is $\sigma$-order continuous, if and only if $E$ has \SDP.
\end{proposition}

\begin{proof}
Normality and the interior point make the norm equivalent to the order unit norm $\|\cdot\|_u$ \cite[Lemma~2.5 and Theorems~2.40, 2.55]{AT07}, and $E$ is Archimedean because the cone is closed. If $\|x_n\|\to 0$, then $c_n=\sup_{m\ge n}\|x_m\|_u\,u$ decreases to $0$ and $\pm x_n\le c_n$. If $\pm x_n\le c_n\downarrow 0$ and the norm is $\sigma$-order continuous, then $\|x_n\|\le(2C+1)\|c_n\|\to 0$; conversely, if order null sequences are norm null, then in particular $\|z_n\|\to 0$ whenever $z_n\downarrow 0$. The last equivalence is Theorem~\ref{thm:banach}.
\end{proof}

Without an order unit, norm null sequences need not be order bounded, e.g.\ $n^{-1/p}e_n$ in $\ell^p$. In $c_0(\Gamma)$ they are order null: if $\|x_n\|\to 0$, then $c_n=\sup_{m\ge n}|x_m|$ exists as the norm limit of finite suprema, decreases and has norm $\sup_{m\ge n}\|x_m\|\to 0$, so $c_n\downarrow 0$. Hence order and norm convergence of sequences coincide in $c_0(\Gamma)$ as well.

Spin factors provide examples without lattice structure in every dimension.

\begin{example}[spin factors]\label{ex:spin}
Let $H$ be a real Hilbert space of dimension at least $2$ and $V=\R\oplus H$ with the Lorentz cone $V_+=\{(a,h):a\ge\|h\|\}$ and the norm $|a|+\|h\|$, which is the order unit norm of $u=(1,0)$ \cite{HOS84}. The cone is closed and normal, $u$ is an interior point, and $V$ is not a vector lattice. The norm is $\sigma$-order continuous. Indeed, let $z_n=(a_n,h_n)\downarrow 0$. Since $z_n-z_{n+1}\in V_+$, $(a_n)$ decreases to some $a\ge 0$ and $\sum_n\|h_n-h_{n+1}\|\le a_1-a$, so $h_n\to h$ and $z_n\to z=(a,h)$ in norm. As $V_+$ is closed, $z\in V_+$ and $z\le z_n$ for all $n$; hence $z=0$. By Proposition~\ref{prop:unit}, every spin factor has \SDP, and order convergence of sequences in it coincides with norm convergence. Since $\dim H$ is arbitrary, this gives non-lattice examples of every infinite dimension.
\end{example}

Operator algebras show both kinds of behaviour.

\begin{example}[operator algebras]\label{ex:vN}
Let $A$ be a unital C$^*$-algebra and $A_{\mathrm{sa}}$ its self-adjoint part with the cone of positive elements, which satisfies the hypotheses of Proposition~\ref{prop:unit} with $u=1$.
\begin{enumerate}[label=(\alph*)]
\item If $A=C(K)$, then $A_{\mathrm{sa}}$ has \SDP\ if and only if $K$ is an almost $P$-space (Theorem~\ref{thm:CK}).
\item If $A$ is an infinite-dimensional von Neumann algebra on $H$, then $A_{\mathrm{sa}}$ fails \SDP. Indeed, $A$ contains a sequence $(p_k)$ of nonzero mutually orthogonal projections. The projections $P_n=\sum_{k\ge n}p_k$ (strong sums) decrease, and $\inf_nP_n=0$ in $A_{\mathrm{sa}}$, because a lower bound $T$ satisfies $\langle T\zeta,\zeta\rangle\le\langle P_n\zeta,\zeta\rangle\to 0$ for every $\zeta\in H$. But $\|P_n\|=1$. This applies to $B(H)_{\mathrm{sa}}$ for infinite-dimensional $H$ and to every infinite-dimensional $L^\infty(\mu)$.
\end{enumerate}
\end{example}

When order convergence and norm convergence of sequences coincide, vector metric spaces are metric spaces in disguise.

\begin{theorem}\label{thm:vm}
Let $E$ be an ordered Banach space with closed normal cone in which order convergence of sequences coincides with norm convergence, and let $\|\cdot\|'$ be an equivalent norm that is monotone on $E_+$; it exists by normality \cite[Theorem~2.38]{AT07}. For an $E$-valued vector metric space $(X,d)$ let $\rho(x,y)=\|d(x,y)\|'$. Then:
\begin{enumerate}[label=\textup{(\roman*)}]
\item $\rho$ is a metric, and $x_n\to x$ if and only if $\rho(x_n,x)\to 0$;
\item $X$ is complete if and only if $(X,\rho)$ is complete;
\item $(X,d)$ has a completion, unique up to isometry; it is the metric completion of $(X,\rho)$ with the continuous extension of $d$.
\end{enumerate}
Conversely, if $E$ is a Banach lattice and \textup{(i)} holds for every $E$-valued vector metric space, then order and norm convergence of sequences coincide in $E$.
\end{theorem}

\begin{proof}
(i) $\rho$ is a metric because $\|\cdot\|'$ is monotone on $E_+$ and $0\le d(x,y)\le d(x,z)+d(z,y)$. Since $d(x_n,x)\ge 0$, $x_n\to x$ means that $(d(x_n,x))$ is order null, which by hypothesis means $\|d(x_n,x)\|'\to 0$.

(ii) A Cauchy sequence for $d$ is $\rho$-Cauchy, since $\|d(x_n,x_{n+p})\|'\le\|a_n\|'\to 0$. A $\rho$-Cauchy sequence has a subsequence with $\sum_k\rho(x_{n_k},x_{n_{k+1}})<\infty$; then $a_k=\sum_{j\ge k}d(x_{n_j},x_{n_{j+1}})$ converges in $E$, decreases to $0$ because $E_+$ is closed and $\|a_k\|\to 0$, and dominates $d(x_{n_k},x_{n_{k+p}})$, so this subsequence is Cauchy for $d$. If $X$ is complete for $d$, the subsequence converges, hence $\rho$-converges by (i), and therefore so does $(x_n)$. If $(X,\rho)$ is complete and $(x_n)$ is Cauchy for $d$, then $\rho(x_n,x)\to 0$ for some $x$, and $x_n\to x$ by (i).

(iii) By hypothesis, order null sequences are norm null, so the norm is $\sigma$-order continuous, and Theorems~\ref{thm:banach} and \ref{thm:main} give a completion, unique up to isometry. Applying (i) and (ii) to it shows that it is a complete metric space for $\|\hat d\|'$ in which $X$ is dense, i.e.\ the metric completion of $(X,\rho)$; $\hat d$ is determined by continuity.

For the converse, take $X=E$ and $d(x,y)=|x-y|$.
\end{proof}

Without lattice structure, one half of the converse still holds. Let $E$ be an ordered Banach space with closed normal cone such that, for every countable $E$-valued vector metric space, convergence of sequences is induced by the metric $\|d\|$. This convergence is then induced by a first countable topology, so $E$ has \SDP\ by Corollary~\ref{cor:nontop}, and its norm is $\sigma$-order continuous by Theorem~\ref{thm:banach}: order null sequences are norm null.

We now turn to spaces of continuous functions on compact spaces. For a compact space $K$, the space $C(K)$ of continuous real functions on $K$, with the pointwise order and the supremum norm, is a Banach lattice with order unit $\mathbf 1$. A Tychonoff space is an \emph{almost $P$-space} if every nonempty $G_\delta$ set, equivalently every nonempty zero set, has nonempty interior \cite{Lev77}.

\begin{theorem}\label{thm:CK}
For a compact Hausdorff space $K$ the following are equivalent.
\begin{enumerate}[label=\textup{(\alph*)}]
\item $K$ is an almost $P$-space.
\item The supremum norm of $C(K)$ is $\sigma$-order continuous.
\item Order convergence of sequences in $C(K)$ is uniform convergence.
\item $C(K)$ has \SDP.
\item Order convergence of sequences in $C(K)$ satisfies the iterated limit axiom.
\item For every $C(K)$-valued vector metric space, $x_n\to x$ if and only if $\|d(x_n,x)\|_\infty\to 0$.
\item For every countable $C(K)$-valued vector metric space, convergence of sequences is induced by a metric.
\item Every countable $C(K)$-valued vector metric space has a completion.
\item Every $C(K)$-valued vector metric space has a completion.
\end{enumerate}
\end{theorem}

The equivalence of (a)--(d) was announced without proof by Veksler \cite[Theorem~2]{Vek73}, in the language of $P'$-spaces; for compact spaces these are exactly the almost $P$-spaces \cite[Theorem~1]{Vek73}. Veksler also lists the condition that every order bounded linear functional on $C(K)$ is $\sigma$-order continuous; compare Theorem~\ref{thm:banach}. A proof of the equivalence of (a), (c) and (d) is given by Hager and van Mill \cite[Theorem~6.4]{HvM15}; see also \cite[Remark~6.5(a)]{HvM15}. We include the short proof of (a)$\Leftrightarrow$(b).

\begin{proof}
(a)$\Rightarrow$(b). Let $x_n\downarrow 0$ with pointwise infimum $g$, and let $c>0$. The set $S_c=\{g\ge c\}=\bigcap_n\{x_n>c-1/n\}$ is a closed $G_\delta$ set, hence a zero set: if $S_c=\bigcap_kU_k$ with $U_k$ open, Urysohn's lemma gives $f_k\in C(K)$ with $0\le f_k\le 1$, $f_k=0$ on $S_c$ and $f_k=1$ off $U_k$, and $S_c$ is the zero set of $\sum_k2^{-k}f_k$. If $S_c$ had nonempty interior $U$, a function $h\in C(K)$ with $0\le h\le c$, $h\ne 0$, vanishing off $U$ would be a positive lower bound of $(x_n)$. So $S_c=\emptyset$ by (a), $g=0$, and $\|x_n\|_\infty\to 0$ by Dini's theorem.

(b)$\Rightarrow$(a). Suppose that (a) fails, and let $Z=f^{-1}(0)$ be a nonempty zero set with empty interior, where $0\le f\le 1$. The functions $\varphi_n=(\mathbf 1-nf)^+$ decrease, and every lower bound $h$ of $(\varphi_n)$ satisfies $h\le 0$ on the dense set $K\setminus Z$, so $h\le 0$ and $\varphi_n\downarrow 0$; but $\|\varphi_n\|_\infty=1$, so (b) fails.

(b)$\Leftrightarrow$(c) is Proposition~\ref{prop:unit} with $u=\mathbf 1$, and (b)$\Leftrightarrow$(d) is Corollary~\ref{cor:BL}. (d)$\Leftrightarrow$(e) is Proposition~\ref{prop:iterated}. (c)$\Rightarrow$(f) is Theorem~\ref{thm:vm}(i), as the supremum norm is monotone; (f)$\Rightarrow$(g) because $\|d\|_\infty$ is a metric; (g)$\Rightarrow$(d) by Corollary~\ref{cor:nontop}. (b)$\Rightarrow$(i) by Theorem~\ref{thm:banach}, (i)$\Rightarrow$(h) trivially, and (h)$\Rightarrow$(d) by Theorem~\ref{thm:main}.
\end{proof}

For all other compact spaces, the negative results of Section~\ref{sec:completions} apply.

\begin{corollary}\label{cor:CKbad}
Let $K$ be a compact space that is not an almost $P$-space, for instance an infinite metrisable space, $\beta\N$, or $[0,\omega_1]$. Then there are countable $C(K)$-valued vector metric spaces without a completion, and countable $C(K)$-valued vector metric spaces whose convergence is induced by no first countable topology.
\end{corollary}

\begin{proof}
The first assertion follows from Theorems~\ref{thm:CK} and \ref{thm:main}, the second from Theorem~\ref{thm:CK} and Corollary~\ref{cor:nontop}. An infinite metrisable $K$ has a non-isolated point $t$, and $\{t\}$ is a zero set with empty interior; in $\beta\N$, $\N^*=\bigcap_n\overline{\{m:m\ge n\}}$ is a zero set with empty interior; in $[0,\omega_1]$, $\{\omega\}=\bigcap_n(n,\omega]$ is a $G_\delta$ set with empty interior.
\end{proof}

An explicit witness: if $Z=f^{-1}(0)$ is a nonempty zero set with empty interior and $\varphi_n=(\mathbf 1-nf)^+$, then the family $a^k_n=k\varphi_n$ has no diagonal, since $k\varphi_{n_k}\le e_k\le e_1$ forces $e_1\ge k$ on $Z$ for all $k$; the spaces of Lemma~\ref{lem:rows} built from it as in Example~\ref{ex:linfty} have the stated properties. On the positive side, compact almost $P$-spaces include the one-point compactifications of uncountable discrete sets (Example~\ref{ex:good}(c)), $\N^*=\beta\N\setminus\N$ \cite{vM84}, and more generally $\beta Y\setminus Y$ for locally compact realcompact $Y$ \cite[\S2]{Vek73}, \cite[Remark~6.5(c)]{HvM15}. Thus $C(\N^*)=\ell^\infty/c_0$ satisfies the conditions of Theorem~\ref{thm:CK}, while $\ell^\infty=C(\beta\N)$ does not.

\section{Products, the bounding number and the spaces \texorpdfstring{$C(X)$}{C(X)}}\label{sec:products}

A product $\prod_{\alpha\in A}E_\alpha$ carries the coordinatewise order; a sequence decreases to $0$ in it if and only if every coordinate does. For $E_\alpha=\R$, Corollary~\ref{cor:RA} below is due to Hager and van Mill \cite[Remark~2.2(ii)]{HvM15}, who derive it from Hager's computation of the $\sigma$-property of $\R^A$ \cite[Theorem~2.1]{Hag16} and their theorem that the $\sigma$-property and \OCRU\ are equivalent for $\R^A$ \cite[Theorem~5.1]{HvM15}. The weaker Egoroff property of $\R^A$ holds if and only if $|A|<\fb$ by \cite{BJ86}; see \cite[\S8]{HvM15} and \cite[Remark~2.5(c)]{Hag16}. Our proofs are direct and use no lattice operations.

\begin{theorem}\label{thm:small}
Let $|A|<\fb$ and let $E_\alpha$ have \SDP\ for each $\alpha\in A$. Then $\prod_{\alpha\in A}E_\alpha$ has \SDP.
\end{theorem}

\begin{proof}
Let $(a^k_n)$ be a family in the product with $a^k_n\downarrow_n 0$, with coordinates $a^k_n(\alpha)$. For each $\alpha$ choose a diagonal $n^\alpha\in\N^{\N}$, $(e^\alpha_k)$ for $(a^k_n(\alpha))$ in $E_\alpha$. Since $|A|<\fb$, there is $g\in\N^{\N}$ with $n^\alpha\le^*g$ for all $\alpha$; the sets $F_\alpha=\{k:g(k)<n^\alpha_k\}$ are finite. Put $n_k=g(k)$ and
\[
e_k(\alpha)=e^\alpha_k+\sum_{j\in F_\alpha,\ j\ge k}a^j_{g(j)}(\alpha).
\]
If $k\notin F_\alpha$, then $a^k_{n_k}(\alpha)\le a^k_{n^\alpha_k}(\alpha)\le e^\alpha_k\le e_k(\alpha)$ by Lemma~\ref{lem:basic}(i). If $k\in F_\alpha$, the term $a^k_{n_k}(\alpha)$ occurs in the sum and the other terms are positive. Each $(e_k(\alpha))_k$ decreases and equals $e^\alpha_k$ for $k>\max F_\alpha$, so $e_k(\alpha)\downarrow 0$. Hence $e_k\downarrow 0$ and $a^k_{n_k}\le e_k$.
\end{proof}

The bound on the number of factors is sharp.

\begin{theorem}\label{thm:large}
Let $|A|\ge\fb$ and suppose that each $E_\alpha$ has a nonzero positive element. Then $\prod_{\alpha\in A}E_\alpha$ fails \SDP.
\end{theorem}

\begin{proof}
Let $\{f_\beta:\beta<\fb\}$ be an unbounded family in $\N^{\N}$, and choose distinct $\alpha_\beta\in A$ and $u_\beta>0$ in $E_{\alpha_\beta}$. Put $a^k_n(\alpha_\beta)=u_\beta$ if $n\le f_\beta(k)$ and $a^k_n(\alpha_\beta)=0$ otherwise, and $a^k_n(\alpha)=0$ for the remaining $\alpha$. Every coordinate of $(a^k_n)_n$ is eventually $0$, so $a^k_n\downarrow_n0$. If $(n_k),(e_k)$ were a diagonal, then, as $(n_k)$ does not bound the family, some $\beta$ satisfies $f_\beta(k)>n_k$ for infinitely many $k$. For these $k$, $e_k(\alpha_\beta)\ge u_\beta$, and since $(e_k)$ decreases this holds for all $k$, contradicting $\inf_ke_k(\alpha_\beta)=0$.
\end{proof}

The construction is a variant of Holbrook's example of a space $\R^X$ without the Egoroff property, with $X=\N^{\N}$ \cite[Example~67.6(v)]{LZ71}.

The hypothesis on the cones cannot be dropped: if the cone of $E$ is $\{0\}$, decreasing sequences are constant, so $E$ has \SDP\ trivially, and so does every product of such spaces.

For powers of the real line we obtain the following.

\begin{corollary}\label{cor:RA}
The space $\R^A$ has \SDP\ if and only if $|A|<\fb$ \cite[Remark~2.2(ii)]{HvM15}. Consequently every $\R^A$-valued vector metric space has a completion if and only if $|A|<\fb$. In particular this holds for $A=\N$.
\end{corollary}

\begin{proof}
This follows from Theorems~\ref{thm:small}, \ref{thm:large} and \ref{thm:main}, since $\R^A$ satisfies \SOC\ by Proposition~\ref{prop:dedekind}.
\end{proof}

Since ZFC does not decide the value of $\fb$, it does not decide the existence of these completions either.

\begin{corollary}\label{cor:indep}
Each of the statements ``$\R^{\omega_1}$ has \SDP'' and ``every $\R^{\omega_1}$-valued vector metric space has a completion'' is independent of ZFC.
\end{corollary}

\begin{proof}
Under CH, $\fb=\aleph_1$ and both statements fail; under MA$+\neg$CH, $\fb=2^{\aleph_0}>\aleph_1$ and both hold \cite{Bla10}. The independence of the first statement also follows from \cite[Remark~2.2(ii)]{HvM15}; under CH even the Egoroff property of $\R^{\omega_1}$ fails \cite{BJ86}, \cite[Theorem~75.3]{LZ71}.
\end{proof}

Spaces of measurable functions over $\sigma$-finite measures behave like countable products.

\begin{proposition}\label{prop:L0}
If $\mu$ is $\sigma$-finite, then $L^0(\mu)$ has \SDP\ and \SOC.
\end{proposition}

This is classical: by \cite[Theorems~71.4 and 71.5]{LZ71}, $L^0(\mu)$ has the $\sigma$-property and \OCRU\ (see also \cite[Remark~5.3(c)]{HvM15}), and the diagonal property is a form of Egoroff's theorem \cite{Hol68}. We include the short direct proof.

\begin{proof}
Property \SOC\ holds by Proposition~\ref{prop:dedekind}. Let $P$ be a probability measure equivalent to $\mu$. For countable families in $L^0(\mu)$ order infima are almost everywhere pointwise infima, so $a^k_n\downarrow_n0$ means $a^k_n\to 0$ almost everywhere, hence in probability. Choose $n_k$ with $P(a^k_{n_k}>2^{-k})<2^{-k}$. By the Borel--Cantelli lemma, almost every point satisfies $a^k_{n_k}\le 2^{-k}$ for large $k$, so $e_k=\sup_{j\ge k}a^j_{n_j}$ is finite almost everywhere, $(e_k)$ decreases, $e_k\to 0$ almost everywhere, and $a^k_{n_k}\le e_k$.
\end{proof}

\begin{remark}\label{rem:L0}
Since $\R^A=L^0(A,\text{counting measure})$, Corollary~\ref{cor:RA} shows that the $\sigma$-finiteness in Proposition~\ref{prop:L0} cannot be dropped in ZFC. On the other hand, it is consistent that $L^0$ of a non-$\sigma$-finite measure has \SDP, namely $\R^{\omega_1}$ when $\aleph_1<\fb$.
\end{remark}

We conclude this section with spaces of continuous functions on general Tychonoff spaces. For a Tychonoff space $X$ we denote by $C(X)$ the vector lattice of all continuous real functions on $X$ with the pointwise order; from here to the end of the section, $X$ is a topological space, not a vector metric space. For compact $X$ this is the Banach lattice $C(K)$ of Theorem~\ref{thm:CK} with $K=X$. In general, $C(X)$ contains unbounded functions and carries no natural norm, so Theorem~\ref{thm:banach} does not apply and the diagonal property has to be decided by order-theoretic means. Its diagonal property has been studied in detail by Hager and van Mill \cite{HvM15} and Hager \cite{Hag16}, in the form of the $\sigma$-property and \OCRU\ (Proposition~\ref{prop:egoroff}). A \emph{$P$-space} is a space in which every $G_\delta$ set is open, and the Lindel\"of number $L(X)$ is the least cardinal $\kappa$ such that every open cover of $X$ has a subcover of cardinality less than $\kappa$. This differs from the usual definition of the Lindel\"of number, which allows subcovers of cardinality at most $\kappa$: with the present convention a Lindel\"of space satisfies $L(X)\le\aleph_1$, and an infinite discrete space satisfies $L(X)=|X|^+$.

\begin{remark}\label{rem:CXknown}
Let $X$ be a Tychonoff space.
\begin{enumerate}[label=(\alph*)]
\item If $C(X)$ has \SDP, then $X$ is an almost $P$-space, by Proposition~\ref{prop:egoroff} and \cite[Theorem~6.3]{HvM15}. If $X$ is an almost $P$-space, then $C(X)$ has \SDP\ if and only if it has the $\sigma$-property \cite[Proposition~6.1 and Theorem~5.1]{HvM15}.
\item If $X$ is a $P$-space with $L(X)\le\fb$, in particular a Lindel\"of $P$-space, then $C(X)$ has the $\sigma$-property \cite[Corollary~3.4(b)]{Hag16}, hence \SDP\ by (a). If $X$ is paracompact and $C(X)$ has the $\sigma$-property, then $L(X)\le\fb$ \cite[Theorem~4.1]{Hag16}. So, for a paracompact $P$-space $X$, the space $C(X)$ has \SDP\ if and only if $L(X)\le\fb$; for discrete $X$ this is Corollary~\ref{cor:RA}.
\item For a topological sum $X$ of compact almost $P$-spaces $X_i$, $i\in I$, the space $C(X)$ has \SDP\ if and only if $|I|<\fb$ \cite[Corollary~9.3]{HvM15}. For $\aleph_1$ copies of $\N^*$ this is independent of ZFC.
\item It is open whether the $\sigma$-property of $C(X)$ implies $L(X)\le\fb$ for realcompact $X$ or for $P$-spaces $X$ \cite[Questions~2.4 and 3.5]{Hag16}, and whether $C(\bigoplus_iX_i)$ has the $\sigma$-property whenever all $C(X_i)$ have it and $|I|<\fb$ \cite[Question~3.2]{Hag16}; see also \cite[Conjecture~9.4]{HvM15}.
\end{enumerate}
\end{remark}

Together with Theorem~\ref{thm:main}, these results give the following.

\begin{corollary}\label{cor:CX}
If $C(X)$ has \SDP, then it satisfies \SOC, and every $C(X)$-valued vector metric space has a completion. This applies if $X$ is a $P$-space with $L(X)\le\fb$, or a topological sum of fewer than $\fb$ compact almost $P$-spaces. Conversely, if $X$ is paracompact and every $C(X)$-valued vector metric space has a completion, then $L(X)\le\fb$.
\end{corollary}

\begin{proof}
The space $C(X)$ is uniformly complete: if $\pm(z_n-z_m)\le\eps_nw$ for $m\ge n$ with $\eps_n\to 0$, then the pointwise limit $z$ satisfies $|z_n-z|\le\eps_nw$, so the convergence is locally uniform and $z$ is continuous. Hence \SOC\ holds by Lemma~\ref{lem:SOCuniform}, and completions exist by Theorem~\ref{thm:main}. The examples are those of Remark~\ref{rem:CXknown}(b),(c). For the converse, Theorem~\ref{thm:main} and Proposition~\ref{prop:egoroff} give the $\sigma$-property, and \cite[Theorem~4.1]{Hag16} gives $L(X)\le\fb$.
\end{proof}

\section{Permanence}\label{sec:perm}

We say that $E$ is the \emph{order direct sum} of subspaces $E_1,E_2$ if $E=E_1\oplus E_2$ and $E_+=(E_+\cap E_1)+(E_+\cap E_2)$, i.e.\ the two projections are positive; projection bands of vector lattices are examples.

\begin{proposition}\label{prop:summand}
If $E=E_1\oplus E_2$ is an order direct sum and $E$ has \SDP, then $E_1$ has \SDP.
\end{proposition}

\begin{proof}
Let $P$ be the projection onto $E_1$. If $a^k_n\downarrow_n0$ in $E_1$, then also in $E$: a lower bound $h\in E$ satisfies $Ph\le Pa^k_n=a^k_n$, so $Ph\le 0$, and $(I-P)h\le(I-P)a^k_n=0$; hence $h\le 0$. Take a diagonal $(n_k),(e_k)$ in $E$. Then $a^k_{n_k}=Pa^k_{n_k}\le Pe_k$, $(Pe_k)$ decreases, and a lower bound $h\in E_1$ of $(Pe_k)$ satisfies $h\le Pe_k\le e_k$ (as $e_k\ge 0$), so $h\le 0$.
\end{proof}

Ideals behave differently from direct summands.

\begin{example}\label{ex:ideals}
The space $c_0$ is an order dense ideal of $\ell^\infty$, and $\ell^\infty$ is an order dense ideal of $\R^{\N}$; here an ideal $I$ of a vector lattice $L$ is \emph{order dense} if for every $y>0$ in $L$ there is $x\in I$ with $0<x\le y$. Since $c_0$ and $\R^{\N}$ have \SDP\ while $\ell^\infty$ does not, \SDP\ passes neither to nor from order dense ideals. For the strong Egoroff property, Luxemburg and Zaanen show with $\ell^\infty$ that it is not inherited by ideals \cite[p.~473]{LZ71}, while $\R^{\N}$ has it \cite[Example~68.6(iv)]{LZ71}.
\end{example}

We now consider pre-Riesz spaces and their vector lattice covers. An ordered vector space $E$ is a \emph{pre-Riesz space} if, for all $x,y,z\in E$, $\{x+z,y+z\}^u\subseteq\{x,y\}^u$ implies $z\in E_+$, where $A^u$ denotes the set of upper bounds of $A$; every directed (that is, $E=E_+-E_+$) Archimedean ordered vector space is pre-Riesz \cite{vH93,KvG19}. A \emph{vector lattice cover} of $E$ is a pair $(Y,i)$ of a vector lattice $Y$ and a bipositive linear map $i\colon E\to Y$ (that is, $x\ge 0$ if and only if $i(x)\ge 0$) such that $y=\inf\{i(x):x\in E,\ y\le i(x)\}$ for every $y\in Y$; by van Haandel's theorem, $E$ has a vector lattice cover if and only if it is pre-Riesz \cite{vH93}, \cite[Section~2.4]{KvG19}. We call $i(E)$ \emph{uniformly dense} in $Y$ if there is $u\in E_+$ such that for all $y\in Y$ and $\eps>0$ there is $x\in E$ with $y\le i(x)\le y+\eps\,i(u)$.

The following lemma collects what we need about covers.

\begin{lemma}\label{lem:cover}
Let $(Y,i)$ be a vector lattice cover of an Archimedean pre-Riesz space $E$.
\begin{enumerate}[label=\textup{(\roman*)}]
\item For a sequence $(c_n)$ in $E$, $c_n\downarrow 0$ in $E$ if and only if $i(c_n)\downarrow 0$ in $Y$.
\item $Y$ is Archimedean.
\item If $i(E)$ is uniformly dense in $Y$ and $b_n\downarrow 0$ in $Y$, then there is a sequence $a_n\downarrow 0$ in $E$ with $b_n\le i(a_n)$ for all $n$.
\end{enumerate}
\end{lemma}

\begin{proof}
Applying the defining property to $-y$ gives $y=\sup\{i(x):i(x)\le y\}$.

(i) Let $c_n\downarrow 0$ in $E$ and let $y$ be a lower bound of $(i(c_n))$. Every $x\in E$ with $i(x)\le y$ is a lower bound of $(c_n)$, so $x\le 0$; hence $y\le 0$. The converse follows from bipositivity.

(ii) Let $ny\le z$ for all $n$, and choose $w\in E$ with $z\le i(w)$, which exists because the set in the definition has infimum $z$ and is therefore nonempty. If $i(x)\le y$, then $nx\le w$ for all $n$, so $x\le 0$; hence $y\le 0$.

(iii) Choose $a_n\in E$ with $b_n+2^{-n}i(u)\le i(a_n)\le b_n+2^{1-n}i(u)$. Then $i(a_{n+1})\le b_{n+1}+2^{-n}i(u)\le b_n+2^{-n}i(u)\le i(a_n)$, so $(a_n)$ decreases. If $h\in E$ is a lower bound, then for $n\ge m$, $i(h)\le b_m+2^{1-n}i(u)$. Taking $n\ge\log_2k$ gives $k(i(h)-b_m)\le 2i(u)$ for all $k$, so $i(h)\le b_m$ by (ii). Hence $i(h)\le 0$ and $h\le 0$.
\end{proof}

Under uniform density, \SDP\ can be tested in the cover.

\begin{theorem}\label{thm:cover}
Let $E$ be an Archimedean pre-Riesz space and $(Y,i)$ a vector lattice cover in which $i(E)$ is uniformly dense. Then $E$ has \SDP\ if and only if $Y$ has \SDP.
\end{theorem}

\begin{proof}
Assume $E$ has \SDP\ and let $a^k_n\downarrow_n0$ in $Y$. By Lemma~\ref{lem:cover}(iii) choose $b^k_n\downarrow_n0$ in $E$ with $a^k_n\le i(b^k_n)$, and a diagonal $(n_k),(e_k)$ for $(b^k_n)$. Then $a^k_{n_k}\le i(e_k)$ and $i(e_k)\downarrow 0$ by Lemma~\ref{lem:cover}(i).

Conversely, let $b^k_n\downarrow_n0$ in $E$. Then $i(b^k_n)\downarrow_n0$ in $Y$; take a diagonal $(n_k),(f_k)$ in $Y$ and, by Lemma~\ref{lem:cover}(iii), $e_k\downarrow 0$ in $E$ with $f_k\le i(e_k)$. Then $b^k_{n_k}\le e_k$ by bipositivity.
\end{proof}

Dense subspaces of $C(K)$ are typical examples.

\begin{example}\label{ex:preriesz}
Let $K$ be a compact Hausdorff space and $E\subseteq C(K)$ a subspace containing $\mathbf 1$ and dense for the supremum norm. Given $f\in C(K)$ and $\eps>0$, density gives $p\in E$ with $f+\eps\le p\le f+3\eps$. Hence $(C(K),\mathrm{incl})$ is a vector lattice cover in which $E$ is uniformly dense, and $E$ has \SDP\ if and only if $C(K)$ has. For example, $C^1[0,1]$ and the polynomials on $[0,1]$, which are pre-Riesz spaces but not vector lattices, fail \SDP\ because $C[0,1]$ does. More generally, by Theorem~\ref{thm:CK}, such a subspace $E$ has \SDP\ if and only if $K$ is an almost $P$-space. For vector sublattices, Hager and van Mill observe that every point-separating vector sublattice of $C(K)$, where $K$ is a compact almost $P$-space, is strongly Egoroff \cite[Remark~6.5(d)]{HvM15}.
\end{example}

Uniform density cannot be replaced by pervasiveness, the other standard density condition for covers.

\begin{example}[pervasiveness does not suffice]\label{ex:pervasive}
Call $i(E)$ \emph{pervasive} in $Y$ if for every $y>0$ in $Y$ there is $x\in E$ with $0<i(x)\le y$. Let $\Gamma$ be uncountable, $Y=\ell^\infty(\Gamma)$, and $E=\R\mathbf 1\oplus c_{00}(\Gamma)$, the functions that are constant off a finite set, with the inclusion. Then $E$ is a vector lattice and $Y$ is a vector lattice cover of $E$: for $y\in Y$ and $t\in\Gamma$, $x_t=\|y\|_\infty\mathbf 1+(y(t)-\|y\|_\infty)e_t\in E$ satisfies $x_t\ge y$ and $x_t(t)=y(t)$. Moreover, $E$ is pervasive, since $0<y(t)e_t\le y$ whenever $y(t)>0$. The space $Y$ fails \SDP, because it contains $\ell^\infty$ as a projection band (Proposition~\ref{prop:summand}, Example~\ref{ex:bad}).

But $E$ has \SDP. First, $x_n\downarrow 0$ in $E$ implies $\|x_n\|_\infty\to 0$. Indeed, $\delta e_t\in E$ for all $t\in\Gamma$ and $\delta>0$, so $x_n\downarrow 0$ in $E$ means $x_n\downarrow 0$ pointwise. Write $x_n=\kappa_n\mathbf 1+g_n$ with $g_n$ finitely supported. Off the countable set $\bigcup_n\operatorname{supp}g_n$ we have $x_n=\kappa_n$, so $\kappa_n\to 0$. Given $\eps>0$, fix $n_0$ with $\kappa_{n_0}<\eps$. For $n\ge n_0$ the sets $T_n=\{t:x_n(t)\ge\eps\}$ are finite and decreasing. If all of them were nonempty, their intersection would contain a point $t$ with $x_n(t)\ge\eps$ for all $n\ge n_0$, which is impossible. Hence $\|x_n\|_\infty<\eps$ for large $n$. Now, given $a^k_n\downarrow_n0$, choose $n_k$ with $\|a^k_{n_k}\|_\infty\le 2^{-k}$; then $(n_k),(e_k)$ with $e_k=2^{1-k}\mathbf 1$ is a diagonal.

So under pervasiveness \SDP\ does not pass from $E$ to $Y$; compare Theorem~\ref{thm:cover}. Since $\mathbf 1\in E$ and $E$ is order dense in the Dedekind complete space $Y$, the space $Y$ is the Dedekind completion of $E$; thus, in ZFC, the diagonal property of a vector lattice need not pass to its Dedekind completion. Under CH, Luxemburg and Zaanen relate the diagonal property of the Dedekind completion to the diagonal property and order separability of the space itself \cite[Exercise~75.8]{LZ71}. A related example for the Egoroff property is given in \cite[Exercise~75.7]{LZ71}: the bounded real functions on $[0,1]$ that differ from their value at $0$ in at most countably many points form a Dedekind $\sigma$-complete space with the Egoroff property whose Dedekind completion lacks it. In accordance with this, $E$ is not order separable (a vector lattice is \emph{order separable} if every subset with a supremum contains a countable subset with the same supremum): $\mathbf 1=\sup\{\mathbf 1_F:F\subseteq\Gamma\text{ finite}\}$ in $E$, but if $C\subseteq\Gamma$ is countable, then the elements $\mathbf 1-\mathbf 1_{\{t\}}$, $t\notin C$, are upper bounds of $\{\mathbf 1_F:F\subseteq C\text{ finite}\}$ that do not dominate $\mathbf 1$. The same cover also shows that sequential order convergence need not pass from $Y$ back to $E$ \cite[Example~8.3]{VC1}. For countable $\Gamma$, the space $E$ itself fails \SDP, since $\mathbf 1_{[n,\infty)}\in E$ decreases to $0$ with norm $1$ and the argument of Example~\ref{ex:linfty} applies.
\end{example}

\section{Further examples and open problems}\label{sec:further}

We begin with value spaces in which already \SOC\ fails.

\begin{example}[spaces failing \SOC]\label{ex:SOC}
The spaces $C[0,1]$, $C^1[0,1]$ and $c$ fail \SOC. Let $\theta$ be a smooth odd nondecreasing function with $\theta(s)=\operatorname{sign}s$ for $|s|\ge\frac12$, and $\psi$ a smooth even function, nonincreasing in $|s|$, with $\psi=1$ on $[-\frac12,\frac12]$ and $\psi=0$ outside $[-1,1]$. Then $z_n(t)=\theta(n(t-\frac12))$ and $c_n(t)=2\psi(n(t-\frac12))$ lie in $C^1[0,1]$, $\pm(z_n-z_m)\le c_n$ for $m\ge n$, and $c_n\downarrow 0$ in $C^1[0,1]$ and in $C[0,1]$. An order limit would be a continuous function agreeing with $\operatorname{sign}(t-\frac12)$ on a dense set, which is impossible. For $c$ take $z_n=\mathbf 1_{2\N\cap[1,2n]}$ and $c_n=\mathbf 1_{[2n,\infty)}$; an order limit would have to be $\mathbf 1_{2\N}\notin c$.
\end{example}

We finally examine the converse of Proposition~\ref{prop:functionals}. By Theorem~\ref{thm:banach}, $\sigma$-order continuity of the continuous positive functionals characterises \SDP\ in ordered Banach spaces with closed normal cone. Outside the Banach setting this fails, already for very simple spaces. For $\R^A$ we use the following classical fact, which we prove for completeness.

\begin{lemma}\label{lem:RAdual}
Let $A$ be a set with $|A|\le 2^{\aleph_0}$. Then every positive linear functional $\varphi$ on $\R^A$ has the form $\varphi(x)=\sum_{j=1}^mc_jx(a_j)$ with $c_j\ge 0$ and $a_j\in A$.
\end{lemma}

\begin{proof}
For $B\subseteq A$ let $\mu(B)=\varphi(\mathbf 1_B)$; this defines a finitely additive measure on the subsets of $A$.

(i) There is no infinite sequence of disjoint sets $B_n$ with $\mu(B_n)>0$. Otherwise the function
\[
x=\sum_n\frac{n}{\mu(B_n)}\mathbf 1_{B_n}\in\R^A
\]
would satisfy $\varphi(x)\ge N$ for all $N$.

(ii) Call $B$ an \emph{atom} if $\mu(B)>0$ and $\mu(C)\in\{0,\mu(B)\}$ for all $C\subseteq B$. Every $B$ with $\mu(B)>0$ contains an atom. Otherwise $B$ splits into two sets of positive measure, neither of which contains an atom, and iterating produces a sequence as excluded in (i). By (i), a maximal disjoint family of atoms is finite, say $A_1,\dots,A_m$, and $R=A\setminus\bigcup_jA_j$ has $\mu(R)=0$, since otherwise $R$ would contain an atom. For each $j$, $U_j=\{B:\mu(B\cap A_j)=\mu(A_j)\}$ is an ultrafilter on $A$, and $\mu(B)=\sum_jc_j[B\in U_j]$ with $c_j=\mu(A_j)$, where $[B\in U_j]$ is $1$ if $B\in U_j$ and $0$ otherwise.

(iii) Each $U_j$ is countably complete. Otherwise there are $D_n\in U_j$ with $D_\infty=\bigcap_nD_n\notin U_j$. The sets $B_n=(D_1\cap\dots\cap D_{n-1})\setminus D_n$, $n\ge 1$, together with $D_\infty$ partition $A$ into sets not in $U_j$. Then $x=\sum_nn\mathbf 1_{B_n}$ satisfies $x\ge m\mathbf 1_{\bigcup_{n\ge m}B_n}$, and $\bigcup_{n\ge m}B_n\in U_j$ as the complement of a finite union of non-members. Hence $\varphi(x)\ge mc_j$ for all $m$, a contradiction. Since $|A|\le 2^{\aleph_0}$ we may regard $A$ as a subset of $\{0,1\}^{\N}$. For each $k$, exactly one of the sets $\{a\in A:a(k)=0\}$ and $\{a\in A:a(k)=1\}$ belongs to $U_j$. The intersection of the chosen sets belongs to $U_j$ by countable completeness and has at most one element, so $U_j$ is generated by a point $a_j$, and $\mu=\sum_jc_j\delta_{a_j}$.

(iv) Let $x\ge 0$ and $F=\{a_1,\dots,a_m\}$. Then $\varphi(x\mathbf 1_F)=\sum_jc_jx(a_j)$. For $y=x\mathbf 1_{A\setminus F}$ let $C_n=\{n-1\le y<n\}\setminus F$, so $\mu(C_n)=0$, and put $z=\sum_nn\mathbf 1_{C_n}\ge y$ and $w=\sum_nn^2\mathbf 1_{C_n}$. For every $k$, $z\le k\mathbf 1_{\bigcup_{n\le k}C_n}+k^{-1}w$, so $\varphi(z)\le\varphi(w)/k$. Hence $\varphi(y)\le\varphi(z)=0$, and linearity gives the claim for all $x$.
\end{proof}

\begin{remark}\label{rem:ulam}
The proof shows more generally that the conclusion holds whenever every countably complete ultrafilter on $A$ is principal. This is the case whenever $|A|$ is smaller than the least measurable cardinal, since the least cardinal carrying a nonprincipal countably complete ultrafilter is measurable \cite[Lemma~10.2]{Jec03}, a result going back to Ulam; we only need $|A|=\fb\le 2^{\aleph_0}$.
\end{remark}

We can now describe to what extent Proposition~\ref{prop:functionals} can be reversed.

\begin{proposition}\label{prop:converse}
\begin{enumerate}[label=\textup{(\alph*)}]
\item In $c_{00}$, the space of finitely supported sequences, every positive linear functional is $\sigma$-order continuous, and there is a strictly positive one; but $c_{00}$ fails \SDP.
\item Let $A$ be a set with $\fb\le|A|\le 2^{\aleph_0}$, e.g.\ $|A|=\fb$. Then $\R^A$ is Dedekind complete and satisfies \SOC, its positive linear functionals separate points and are all $\sigma$-order continuous, but $\R^A$ fails \SDP; hence some $\R^A$-valued vector metric space has no completion.
\item In an ordered Banach space with closed normal cone, \SDP\ holds if and only if every positive linear functional is $\sigma$-order continuous, if and only if every continuous positive linear functional is $\sigma$-order continuous.
\end{enumerate}
\end{proposition}

\begin{proof}
(a) Every linear functional has the form $\varphi(x)=\sum_nc_nx(n)$ with $c_n=\varphi(e_n)$, which is $\ge 0$ if $\varphi$ is positive. If $x_k\downarrow 0$, then all $x_k$ are supported in the finite support of $x_1$, so $x_k\to 0$ coordinatewise there and $\varphi(x_k)\to 0$. Take $c_n=2^{-n}$ for a strictly positive one. For $a^k_n=e_k/n$, a diagonal $(n_k),(d_k)$ would satisfy $d_1\ge d_k\ge e_k/n_k$ for all $k$, so $d_1$ would have infinite support. (Alternatively, $c_{00}$ fails the $\sigma$-property, and Proposition~\ref{prop:egoroff} applies; cf.\ \cite[Remark~5.2(i)]{HvM15}.)

(b) Dedekind completeness and \SOC\ are clear (Proposition~\ref{prop:dedekind}). By Lemma~\ref{lem:RAdual}, positive functionals are finite combinations of point evaluations, which separate points and are $\sigma$-order continuous because $x_n\downarrow 0$ in $\R^A$ means pointwise convergence to $0$. Property \SDP\ fails by Theorem~\ref{thm:large}, and completions fail to exist by Theorem~\ref{thm:main}.

(c) Property \SDP\ implies that every positive functional is $\sigma$-order continuous (Proposition~\ref{prop:functionals}); this includes the continuous ones; and $\sigma$-order continuity of the continuous positive functionals implies \SDP\ by Theorem~\ref{thm:banach}. No automatic continuity is needed.
\end{proof}

Thus the converse of Proposition~\ref{prop:functionals} fails even in Dedekind complete spaces with \SOC\ and a separating family of positive functionals. For $\R^A$, the obstruction is purely combinatorial: positive functionals see only finitely many coordinates, whereas \SDP\ requires a single element dominating a diagonal over all coordinates at once.

\begin{table}[ht]
\centering
\caption{Value spaces $E$ and completions of $E$-valued vector metric spaces}\label{tab:summary}
\begin{tabular}{@{}lccc@{}}
\toprule
$E$ & \SOC & \SDP & completions\\
\midrule
finite-dimensional, closed cone & yes & yes & yes\\
$L^p(\mu)$ $(1\le p<\infty)$, $c_0(\Gamma)$, $S_2^{\mathrm{sa}}$, spin factors & yes & yes & yes\\
$C(K)$, $K$ compact almost $P$, e.g.\ $\ell^\infty/c_0$ & yes & yes & yes\\
$C(X)$, $X$ a $P$-space with $L(X)\le\fb$ & yes & yes & yes\\
$\R^A$ with $|A|<\fb$; $L^0(\mu)$, $\mu$ $\sigma$-finite & yes & yes & yes\\
$\R^A$ with $|A|\ge\fb$ & yes & no & no\\
$\ell^\infty$ & yes & no & no\\
$c_{00}$ & yes & no & no\\
$c$, $C[0,1]$, $C^1[0,1]$ & no & no & no\\
\bottomrule
\end{tabular}
\end{table}

\begin{remark}\label{rem:forthcoming}
Spaces of continuous functions are treated in Theorem~\ref{thm:CK} and Corollary~\ref{cor:CX}; for $C(X)$ several questions remain open; see Remark~\ref{rem:CXknown}(d). The analogue of \SDP\ for order convergence of nets and the corresponding completions by Cauchy nets will be treated in forthcoming work.
\end{remark}

\begin{problem}
In vector lattices, \SDP\ is the conjunction of the Egoroff property and the dominating property \cite[Theorem~3.4]{Dod74}, and it is characterised by the pseudo order closures of convex sets \cite[Theorem~3.14]{Dod74}. Find an Egoroff property for ordered vector spaces, for instance in terms of order intervals, for which the first statement remains true, and decide whether the second holds in pre-Riesz spaces. More generally, how is \SDP\ of a pre-Riesz space related to \SDP\ of its Riesz completion, the vector lattice cover generated by the image of the space as a vector lattice \cite{vH93}, when the space is not uniformly dense in it?
\end{problem}

\begin{problem}
By Example~\ref{ex:pervasive}, pervasiveness does not suffice for the passage of \SDP\ from $E$ to $Y$ in Theorem~\ref{thm:cover}; does it suffice for the passage from $Y$ to $E$? For which $E$ is there an order dense extension $F\supseteq E$ with \SOC\ and \SDP, so that all $E$-valued vector metric spaces have $F$-valued completions? For a vector lattice and its Dedekind completion, compare \cite[Exercise~75.8]{LZ71} and Example~\ref{ex:pervasive}.
\end{problem}


\begin{thebibliography}{99}
\bibitem{AT07} C.\,D. Aliprantis, R. Tourky, \emph{Cones and Duality}, Graduate Studies in Mathematics 84, American Mathematical Society, Providence, 2007.
\bibitem{Bla10} A. Blass, Combinatorial cardinal characteristics of the continuum, in: \emph{Handbook of Set Theory}, Springer, Dordrecht, 2010, 395--489.
\bibitem{BJ86} A. Blass, T. Jech, On the Egoroff property of pointwise convergent sequences of functions, \emph{Proc. Amer. Math. Soc.} \textbf{98} (1986), no.~3, 524--526.
\bibitem{Cev14} C. \c{C}evik, On continuity of functions between vector metric spaces, \emph{J. Funct. Spaces} \textbf{2014} (2014), Article ID 753969, 6~pp.
\bibitem{Cev15} C. \c{C}evik, Completion of vector metric spaces, \emph{Gazi Univ. J. Sci.} \textbf{28} (2015), 69--73.
\bibitem{VC1} C. \c{C}evik, Vector metrics with values in pre-Riesz spaces I: transfer of order convergence along vector lattice covers, preprint, arXiv:2611.12376 (2026).
\bibitem{CA09} C. \c{C}evik, I. Altun, Vector metric spaces and some properties, \emph{Topol. Methods Nonlinear Anal.} \textbf{34} (2009), 375--382.
\bibitem{Dod74} T.\,K.\,Y. Chow Dodds, Egoroff properties and the order topology in Riesz spaces, \emph{Trans. Amer. Math. Soc.} \textbf{187} (1974), 365--375.
\bibitem{Hag16} A.\,W. Hager, The $\sigma$-property in $C(X)$, \emph{Comment. Math. Univ. Carolin.} \textbf{57} (2016), 231--239.
\bibitem{HvM15} A.\,W. Hager, J. van Mill, Egoroff, $\sigma$, and convergence properties in some archimedean vector lattices, \emph{Studia Math.} \textbf{231} (2015), 269--285.
\bibitem{HOS84} H. Hanche-Olsen, E. St\o rmer, \emph{Jordan Operator Algebras}, Pitman, Boston, 1984.
\bibitem{Hol68} J.\,A.\,R. Holbrook, Seminorms and the Egoroff property in Riesz spaces, \emph{Trans. Amer. Math. Soc.} \textbf{132} (1968), 67--77.
\bibitem{Jec03} T. Jech, \emph{Set Theory}, third millennium edition, Springer, Berlin, 2003.
\bibitem{Kad51} R.\,V. Kadison, Order properties of bounded self-adjoint operators, \emph{Proc. Amer. Math. Soc.} \textbf{2} (1951), 505--510.
\bibitem{KvG19} A. Kalauch, O. van Gaans, \emph{Pre-Riesz Spaces}, De Gruyter Expositions in Mathematics 66, De Gruyter, Berlin, 2019.
\bibitem{Kan37} L.\,V. Kantorovich, Lineare halbgeordnete R\"aume, \emph{Mat. Sb. (N.S.)} \textbf{2(44)} (1937), no.~1, 121--168.
\bibitem{Lev77} R. Levy, Almost-$P$-spaces, \emph{Canad. J. Math.} \textbf{29} (1977), 284--288.
\bibitem{LZ71} W.\,A.\,J. Luxemburg, A.\,C. Zaanen, \emph{Riesz Spaces I}, North-Holland, Amsterdam, 1971.
\bibitem{OC20} \c{C}.\,C. \"Ozeken, C. \c{C}evik, Unbounded vectorial Cauchy completion of vector metric spaces, \emph{Gazi Univ. J. Sci.} \textbf{33} (2020), 761--765.
\bibitem{vD84} E.\,K. van Douwen, The integers and topology, in: \emph{Handbook of Set-Theoretic Topology}, North-Holland, Amsterdam, 1984, 111--167.
\bibitem{vH93} M. van Haandel, \emph{Completions in Riesz Space Theory}, PhD thesis, University of Nijmegen, 1993.
\bibitem{vM84} J. van Mill, An introduction to $\beta\omega$, in: \emph{Handbook of Set-Theoretic Topology}, North-Holland, Amsterdam, 1984, 503--567.
\bibitem{Vek73} A.\,I. Veksler, $P'$-points, $P'$-sets, $P'$-spaces. A new class of order-continuous measures and functionals, \emph{Dokl. Akad. Nauk SSSR} \textbf{212} (1973), no.~4, 789--792 (in Russian); English transl.: \emph{Soviet Math. Dokl.} \textbf{14} (1973), 1445--1450.
\bibitem{Zaa68} A.\,C. Zaanen, Stability of order convergence and regularity in Riesz spaces, \emph{Studia Math.} \textbf{31} (1968), 159--172.
\bibitem{Zab97} P.\,P. Zabrejko, $K$-metric and $K$-normed linear spaces: survey, \emph{Collect. Math.} \textbf{48} (1997), 825--859.
\end{thebibliography}
\end{document}